\documentclass[12pt]{article}
\usepackage[utf8]{inputenc}
\usepackage[british]{babel}
\usepackage{cmap}
\usepackage{lmodern}

\usepackage{amssymb, amsmath, amsthm}
\usepackage[a4paper,top=25mm,bottom=25mm,left=25mm,right=25mm]{geometry}
\usepackage{ragged2e}

\usepackage{authblk} 
\usepackage{pifont}
\usepackage{graphicx}
\usepackage[dvipsnames,svgnames,table]{xcolor}
\usepackage[figuresright]{rotating}
\usepackage{xtab} 
\usepackage{longtable} 
\usepackage{multirow}
\usepackage{footnote}
\usepackage[stable]{footmisc}
\usepackage{chngpage} 
\usepackage{pdflscape} 
\usepackage[nottoc,notlot,notlof]{tocbibind} 

\usepackage{pgfplots}
\pgfplotsset{every tick label/.append style={font=\footnotesize}}
\pgfplotsset{compat=1.18}
\usepackage{setspace}

\usepackage{array}
\newcolumntype{K}[1]{>{\centering\arraybackslash$}p{#1}<{$}}

\makesavenoteenv{tabular}
\usepackage{tabularx}
\usepackage{booktabs}
\usepackage{threeparttable}
\usepackage[referable]{threeparttablex} 
\newcolumntype{R}{>{\raggedleft\arraybackslash}X}
\newcolumntype{L}{>{\raggedright\arraybackslash}X}
\newcolumntype{C}{>{\centering\arraybackslash}X}
\newcolumntype{A}{>{\columncolor{gray!25}}C}
\newcolumntype{a}{>{\columncolor{gray!25}}c}

\newlength{\tablen}

\usepackage{dcolumn} 
\newcolumntype{.}{D{.}{.}{-1}}

\usepackage{tikz}
\usetikzlibrary{arrows, calc, matrix, patterns, positioning, trees}
\usepackage[semicolon]{natbib} 
\usepackage[hyphens]{xurl}
\usepackage{microtype}
\usepackage[justification=centering]{caption} 

\usepackage[labelformat=simple]{subcaption}

\DeclareCaptionLabelFormat{parenthesis}{(#2)}
\makeatletter
\renewcommand\p@subfigure{\arabic{figure}.}
\makeatother

\DeclareCaptionLabelFormat{parenthesis}{(#2)}
\makeatletter
\renewcommand\p@subtable{\arabic{table}.}
\makeatother

\usepackage{enumitem}

\setlist[itemize]{leftmargin=2.5\parindent}
\setlist[enumerate]{leftmargin=2.5\parindent}

\usepackage{hyperref} 
\hypersetup{
  colorlinks   = true,    		
  urlcolor     = blue,    		
  linkcolor    = blue,    		
  citecolor    = ForestGreen	
}

\def\addlegendimage{\csname pgfplots@addlegendimage\endcsname}

\theoremstyle{plain}

\theoremstyle{definition}

\theoremstyle{remark}

\makeatletter
\let\@fnsymbol\@alph
\makeatother

\def\keywords{\vspace{.5em} 
{\noindent \textit{Keywords}: }}

\def\AMS{\vspace{.5em} 
{\noindent \textbf{\emph{MSC} class}: }}

\def\JEL{\vspace{.5em} 
{\noindent \textbf{\emph{JEL} classification number}: }}

\title{Increasing competitiveness by imbalanced groups: \\ The example of the 48-team FIFA World Cup}

\author{\href{https://sites.google.com/view/laszlocsato}{L\'aszl\'o Csat\'o}\thanks{~Corresponding author. Email: \emph{laszlo.csato@sztaki.hun-ren.hu} \newline
Institute for Computer Science and Control (SZTAKI), Hungarian Research Network (HUN-REN), Laboratory on Engineering and Management Intelligence, Research Group of Operations Research and Decision Systems, Budapest, Hungary \newline
Corvinus University of Budapest (BCE), Institute of Operations and Decision Sciences, Department of Operations Research and Actuarial Sciences, Budapest, Hungary}
$\qquad \qquad$
Andr\'as Gyimesi\thanks{~Email: \emph{gyimesi.andras@ktk.pte.hu} \newline
University of P\'ecs, P\'ecs, Hungary \newline
Institute for Computer Science and Control (SZTAKI), Hungarian Research Network (HUN-REN), Laboratory on Engineering and Management Intelligence, Research Group of Operations Research and Decision Systems, Budapest, Hungary}}

\date{\today}

\def\Dedication{
{\noindent
``\emph{The objectives of FIFA are: [\dots] \\
g) to promote integrity, ethics and fair play with a view to preventing
all methods or practices, such as corruption, doping or match manipulation, which might jeopardise the integrity of matches, competitions, players, officials and member associations or give rise to abuse of association football.}''
}
\vspace{0.25cm}

\flushright
\noindent (FIFA Statutes, May 2022 edition \citep[Article~2]{FIFA2022d})

\vspace{1cm} 
\justify }

\begin{document}

\maketitle

\thispagestyle{empty}
\Dedication

\begin{abstract}
\noindent
A match played in a sports tournament can be called stakeless if at least one team is indifferent to its outcome because it already has qualified or has been eliminated. Such a game threatens fairness since teams may not exert full effort without incentives. This paper suggests a novel classification for stakeless matches based on their expected outcome: they are more costly if the indifferent team is more likely to win by playing honestly.
Our approach is illustrated with the 2026 FIFA World Cup, the first edition of the competition with 48 teams. We propose a novel format based on imbalanced groups, which substantially reduces the probability of stakeless matches played by the strongest teams according to Monte Carlo simulations. The new design also increases the uncertainty of match outcomes and requires fewer matches.
Governing bodies in sports are encouraged to consider our innovative idea in order to enhance the competitiveness of their tournaments.


\AMS{90-10, 90B90, 91B14}

\JEL{C44, D71, Z20}
\end{abstract}

\clearpage

\section{Introduction} \label{Sec1}

On 14 March 2023, the FIFA (F\'ed\'eration Internationale de Football Association, International Federation of Association Football) Council unanimously amended the format of the 2026 FIFA World Cup with the following commentary \citep{FIFA2023b}:
``\emph{The revised format mitigates the risk of collusion and ensures that all the teams play a minimum of three matches, while providing balanced rest time between competing teams.}''
However, collusion is not the only threat to the integrity of sports. The recent academic literature has argued that the absence of win incentives is also unfair \citep{ChaterArrondelGayantLaslier2021, CsatoMolontayPinter2024, DevriesereCsatoGoossens2025, Guyon2022a}. This can be a problem especially in the last round of group matches. For example, both Brazil and Nigeria in the 1998 FIFA World Cup, as well as Italy in the 2016 UEFA European Championship, secured the first place in their group after two rounds---and they all lost their last games \citep{Guyon2022a}.


A match is stakeless for a team if either it has reached its objective (e.g.\ qualification for the next round), or its advance to the next round has become impossible. An illustration of the former case is France in Group D of the 2022 FIFA World Cup. The French coach \emph{Didier Deschamps}, whose team had already qualified for the knockout stage, made nine changes compared to the previous match, and the second squad performed as if they had never played together before \citep{LeMonde2022}. Unsurprisingly, they lost to Tunisia despite being the defending champion. Since Canada was also eliminated in Group F before the last round, both France and Canada played a weakly stakeless match \citep{CsatoMolontayPinter2024} in the sense that their qualification was decided, but this condition did not hold for their opponents, Tunisia and Morocco, respectively.

However, the two matches were fundamentally different: according to the World Football Elo Ratings, a widely used measure of strength in football \citep{ChaterArrondelGayantLaslier2021, Csato2022a, GasquezRoyuela2016, LasekSzlavikBhulai2013}, France was more likely to win than Canada assuming that all teams exert full effort. In particular, the win expectancy of France was 88.35\% (implied by having 352 Elo points more than Tunisia), while the win expectancy of Canada was only 33.51\% (implied by having 119 Elo points less than Morocco). Therefore, the first stakeless match between France and Tunisia can be considered more costly for the organiser as the result would be more surprising and would have a higher impact on the other teams if the indifferent team does not play honestly. Although the actual winning probability of France could have been lower due to the unusual French squad, it had probably been chosen by the coach exactly because the match was stakeless for them.

The existing literature \citep{ChaterArrondelGayantLaslier2021, CsatoMolontayPinter2024} assumes a uniform cost for stakeless matches. One of our main contributions resides in introducing a weighting scheme for non-competitive games.
This idea will be used to demonstrate how the design of a 48-team FIFA World Cup can be improved by creating deliberately imbalanced groups during the draw.

The new format is found to have three crucial advantages, which are detailed in Section~\ref{Sec42}:
\begin{itemize}
\item 
It contains fewer matches, with an expected reduction of 8--9\% in the number of games played by the strongest teams whose players have the highest workload (Figure~\ref{Fig3});
\item 
It implies more uncertain group games for the strongest teams (Figure~\ref{Fig4}), more uncertain games in the knockout stage from the quarterfinals (Figure~\ref{Fig5}), and a higher proportion of games between the strongest teams (Figure~\ref{Fig6});
\item 
It substantially reduces the ratio of stakeless matches for the 16 strongest teams (Figure~\ref{Fig7}).
\end{itemize}
These benefits present powerful arguments for both researchers and tournament organisers to consider imbalanced competition designs.

Our study is structured as follows.
Section~\ref{Sec2} gives a literature overview. Section~\ref{Sec3} describes two alternative formats for the 2026 FIFA World Cup: the official design and our proposal containing imbalanced groups. The simulation technique, the tournament metrics, and the limitations of our study are also discussed in Section~\ref{Sec3}. The results are presented and evaluated in Section~\ref{Sec4}. Finally, Section~\ref{Sec5} provides concluding remarks.

\section{Related literature} \label{Sec2}

Tournament design has received serious attention in both economics \citep{Medcalfe2024, Palacios-Huerta2023b, Szymanski2003} and operational research \citep{Wright2009, Wright2014, KendallLenten2017, Csato2021a, DevriesereCsatoGoossens2025}. 
In particular, the format of the 2026 FIFA World Cup has been extensively investigated. FIFA has approved the expansion of the 2026 FIFA World Cup to 48 teams in January 2017 and decided that the competition would start with 16 groups of three teams each (80 matches in total) \citep{FIFA2017d}. According to \citet{Truta2018}, this format would not have increased the number of non-competitive matches compared to the traditional design with 32 teams, eight groups of four teams each. Nonetheless, a predetermined schedule where the \emph{a priori} strongest team plays against the second strongest team in the first round and the weakest team in the second round can roughly halve the probability of non-competitive matches.

\citet{Guyon2020a} examines the risk of collusion in the same setting: since the two teams playing the last group game know exactly what results will let them qualify for the knockout stage, a particular result can be easily beneficial for both of them at the expense of the third team. The best solution to minimise the risk of collusion turns out to be if the strongest team plays the first two group matches. However, then only the best team is threatened by collusion, and none of the 16 favourites play in the last round of the group stage. Furthermore, the risk of collusion in at least one group becomes unacceptably high. The problem cannot be mitigated by forbidding draws during the group stage or using the 3-2-1-0 point system, which gives an additional point to the winner of the penalty shootout after a draw.
Hence, \citet{Guyon2020a} has proposed seven alternative formats with groups of three, four, or six teams, but has not considered imbalanced groups.
The Monte Carlo simulations of \citet{ChaterArrondelGayantLaslier2021} generally reinforce the results of \citet{Guyon2020a} regarding both the high risk of non-competitive matches and the optimal schedule.

Probably, these findings have inspired FIFA to revise the format of the 2026 FIFA World Cup and decide for 12 groups of four teams each---Format 1 in \citet{Guyon2020a}---in 2023 \citep{FIFA2023b} as mentioned in the Introduction.

Some further designs have been proposed for the 48-team tournament. \citet{Renno-Costa2023} presents a double-elimination structure instead of the group stage in order to produce more competitive and exciting matches such that teams with an early loss can still recover. \citet{GuajardoKrumer2024} recommend three novel tournament formats. The first, inspired by beach volleyball, contains only two rounds of dynamically scheduled group matches (that is, the set of games played in the second round depends on the outcome of the first round). The second is based on the usual design, but the knockout bracket is organised with 24 teams to select three winners, the best directly qualifying for the final, and the other two playing a semifinal. The third proposal is a hybrid of the first two options. The total number of games ranges from 71 to 95, which is considerably lower than the undesirably high 104 games of the current format.

In the following, another format will be suggested for a 48-team FIFA World Cup with imbalanced groups and compared to the official format with 12 groups of four teams each. In particular, we will focus on the distribution of stakeless matches across the teams.

The expansion to 48 teams raises the issue of how these slots are allocated among the FIFA confederations, too. \citet{KrumerMoreno-Ternero2023} use standard tools from the fair allocation literature for this purpose and find that UEFA has a solid basis for claiming additional slots. \citet{CsatoKissSzadoczki2025} extend the methodology of the FIFA World Ranking to evaluate the performance of sets of teams and suggest a transparent slot allocation policy. According to their results, more European and South American teams should play in the (2026) FIFA World Cup.

There are several studies on group balance, a key property for the group draw, in the FIFA World Cup.
According to \citet{Guyon2015a}, the 2014 FIFA World Cup draw has implied a substantial level of group imbalance, the eight groups have not been of the same competitive level. The 2022 FIFA World Cup draw has not guaranteed balance to the extent possible because of the misaligned treatment of the winners of the play-offs \citep{Csato2022d}. Various procedures have been proposed to create balanced groups by taking geographic constraints into account \citep{CeaDuranGuajardoSureSiebertZamorano2020, LalienaLopez2019, LalienaLopez2025}. Group balance is not merely a theoretical curiosity as demonstrated by \citet{LaprePalazzolo2023}: in the FIFA World Cups between 1954 and 2022, the range in the rating of group opponents exceeds 118 Elo points, which is roughly the difference between an average participant and an average semifinalist. Furthermore, an increase of merely 88 Elo points in the rating of group opponents is expected to halve the probability of reaching the quarterfinals in a 32-team tournament.

Our paper is strongly connected to another line of literature that quantifies the probability of collusion and match-fixing.
\citet{Guyon2020a} calculates the risk of collusion in groups of three if two teams advance to the next phase.
\citet{ChaterArrondelGayantLaslier2021} distinguish three types of games: competitive (when neither team is indifferent and their targets are incompatible), collusive (when the targets of both teams are compatible and neither is indifferent), and stakeless (when at least one team is completely indifferent). Their probabilities are determined in various settings of the FIFA World Cup: 8 groups of four teams, 16 groups of three teams, 12 groups of four teams, 8 groups of five teams, and 12 groups of five teams.
\citet{Stronka2024} offers two innovative changes for groups of three teams: random tie-breaking based on goal difference and dynamic scheduling. These proposals are able to reduce the expected number of matches with a high risk of collusion from 5.5 to 0.26 in the 2026 FIFA World Cup.

\citet{CsatoMolontayPinter2024} consider another classification scheme, where the games can be competitive, weakly stakeless (exactly one team is indifferent), and strongly stakeless (both teams are indifferent). Their probabilities are computed for the 12 valid schedules in the UEFA Champions League group stage. Some schedules turn out to be dominated by other schedules. In the optimal schedule, the strongest team should play at home in the last round against one of the middle teams, depending on the preferences of the organiser.
\citet{Gyimesi2024} follows this approach with a particular weighting system (0 for competitive, 0.5 for weakly stakeless, 1 for strongly stakeless) to measure the ratio of stakeless games in the league phase of the UEFA Champions League that has been introduced from the 2024/25 season. The results suggest that the novel design substantially reduces the probability of stakeless games.
\citet{DevriesereGoossensSpieksma2025} identify three distinct types of non-competitive matches: asymmetric (exactly one team is indifferent), stakeless (both teams are indifferent), and collusive (none of the teams is indifferent, and a particular result secures the desired prize for both of them). Their probabilities are quantified in the pre-2024 standard group stage and the post-2024 incomplete round-robin league phase of the UEFA Champions League. The new format is found to perform substantially better under several reasonable schedules.

\citet{Csato2022a} presents a method to quantify the violation of strategy-proofness through the example of the European Qualifiers for the 2022 FIFA World Cup. The threat of tanking can be substantially mitigated by adding a carefully chosen set of constraints to the group draw.
\citet{Csato2023a} demonstrates that tie-breaking rules might affect the occurrence of a situation when the final position of a team is already secured, independently of the results in the last round of group matches.
In particular, merely the tie-breaking policy of the 2024 UEFA European Football Championship has increased the risk of match-fixing by more than 10 percentage points \citep{Csato2025d}. 

Nevertheless, none of the studies above take the identity of the teams that play these non-competitive matches into account. However, as has been argued in the Introduction, stakeless matches are more threatening for the organiser if they affect the best teams: a weak team is more likely to lose the match anyway, thus, the cost of misaligned incentives is smaller. Introducing that aspect is one of our main contributions to the literature.

\section{Methodology} \label{Sec3}

This section summarises the methodology of the study. Section~\ref{Sec31} provides an overview of two tournament formats for the 2026 FIFA World Cup. The simulation framework is detailed in Section~\ref{Sec32}. Section~\ref{Sec33} discusses the underlying data, and Section~\ref{Sec34} defines the metrics used to compare the two alternative designs. Last but not least, the limitations of our results are outlined in Section~\ref{Sec35}.

\subsection{Tournament formats} \label{Sec31}

The official design of the 2026 FIFA World Cup is detailed in Section~\ref{Sec311}, while a novel proposal based on imbalanced groups is specified in Section~\ref{Sec312}. Note that the latter has not appeared in the previous literature, even though \citet{Guyon2020a} studies eight, \citet{Renno-Costa2023} one, and \citet{GuajardoKrumer2024} three reasonable formats.

\subsubsection{The official format} \label{Sec311}

As discussed in Section~\ref{Sec2}, FIFA finally decided to organise the 2026 FIFA World Cup with 12 groups of four teams each, followed by a knockout stage. The first two teams from each group and the eight best third-placed teams qualify for the Round of 32. This structure raises several fairness issues since the groups are \emph{ex post} (after the draw) treated separately; for example, some group winners play against a runner-up and some against a third-placed team in the Round of 32 \citep{Guyon2018a, Csato2021a}. Fortunately, these rules are irrelevant to our simulations since the groups are \emph{ex ante} (before the draw) equivalent.

The groups are played in a single round-robin format. The ranking in each group is determined by the number of points, followed by head-to-head records (head-to-head number of points, goal difference, number of goals scored). The next tie-breaking criteria are goal difference, the number of goals scored, disciplinary points, and position in the latest FIFA Men's World Ranking \citep[Article~13]{FIFA2025}. 

At the moment, it is not known how the seeding will work in the 2026 FIFA World Cup. It is assumed to follow the previous editions \citep{Csato2023d}: the 48 teams are partitioned into four pots of 12 teams each according to the official FIFA World Ranking of November 2025 except for the three host nations (Canada, Mexico, United States), which are automatically assigned to the strongest Pot 1, and for the two winners of the play-offs, which are assigned to the weakest Pot 4. The draw chooses one team from each pot for each group. Two teams from the same confederation (other than UEFA) cannot play in the same group, and the number of European teams should be between one and two in all groups. The two winners of the play-offs are not taken for any confederation into account. Similar to \citet{Stronka2024}, the draw procedure is implemented by a rejection sampler that is---in contrast to the official Skip mechanism \citep{Csato2025c}---uniformly distributed over the set of valid assignments \citep{RobertsRosenthal2024}.

In the knockout stage, the bracket of the 2024 UEFA European Championship \citep{Guyon2018a} is essentially ``doubled'': group G corresponds to Group A, Group H to Group B, and so on. For instance, eight group winners play against the eight third-placed teams, four group winners play against four runners-up, and eight runners-up play against each other in the Round of 32. In particular, the eight best third-placed teams play against the winners of Groups B, C, E, F, H, I, K, L after a random matching. From the Round of 16, the knockout bracket is deterministic.

\subsubsection{Our recommendation: imbalanced groups} \label{Sec312}

The format of the 2026 FIFA World Cup has two substantial shortcomings. First, since the design is analogous to the recent UEFA European Championships, the same fairness issues arise: it is impossible to treat the groups in the knockout stage equally \citep{Guyon2018a, Csato2021a}. Furthermore, both within-group differences (as the number of teams is expanded from 32 to 48) and the ratio of advancing teams (from 1/2 to 2/3) increase compared to the previous editions. These changes are expected to make qualification easier for the stronger teams, implying fewer competitive matches in the group stage.

Therefore, the main motivation behind our proposal is to improve the competitiveness of the group matches. In particular, two types of groups are introduced: eight Tier 1 groups of four stronger teams and four Tier 2 groups of four weaker teams. Instead of the usual Round of 32, an intermediate play-off round is organised, similar to the novel format of UEFA club competitions \citep{Gyimesi2024}. From Tier 1 groups, the group winner directly qualifies for the Round of 16, and the runner-up enters the play-off round. On the other hand, from Tier 2 groups, both the group winner and the runner-up qualify for the play-off round, which is contested by the eight runners-up of Tier 1 groups besides them.

In this play-off round, a Tier 1 runner-up faces either a Tier 2 group winner or a runner-up. Again, a random matching is assumed. Naturally, the winners of the play-off round are paired against the Tier 1 group winners in the Round of 16. From this point on, the format follows the official FIFA World Cup format.

Seeding is similar to the official format with eight pots, four (Pots 1, 2, 5, 7) containing 8 teams and another four (Pots 3, 4, 6, 8) containing 4 teams each. The three hosts are assigned to Pot 1 together with the five strongest teams. The two winners of the play-offs are assigned to Pot 8 together with the two weakest directly qualified teams.
Tier 1 groups receive one team from Pots 1, 2, 5, 7 each, while Tier 2 groups receive one team from Pots 3, 4, 6, 8 each. Following the official format, no group can contain more than one team from the confederations AFC, CAF, CONCACAF, CONMEBOL, OFC, and should contain at least one and at most two UEFA teams, which is guaranteed by a rejection sampler. Again, the two winners of the play-offs are not taken for any confederation into account.

Imbalanced groups have already been used in some sports competitions \citep[Section~5.3]{DevriesereCsatoGoossens2025}. The first example is the EHF (Men's Handball) Champions League between the 2015/16 and 2019/20 seasons, which has been extensively studied by \citet{Csato2020b}. In particular, the design is found to considerably increase the proportion of high quality and more balanced games. UEFA follows a similar principle in its Nations League launched in 2018, where the national teams are divided into four leagues of different strengths \citep{Csato2022a, ScellesFrancoisValenti2024}. Finally, the 2024 European Water Polo Championships have also consisted of four groups in two divisions \citep{LEN2023}, and the 2025 IIHF (Ice Hockey) Women's World Championship has contained two tiered groups \citep{IIHF2024}.

Unfortunately, tournament organisers usually explain their decisions for imbalanced groups only by straightforward arguments such as guaranteeing more matches between the leading teams or increasing attractiveness \citep{EHF2014a}. However, the recent examples in water polo and women's ice hockey construct the imbalanced groups such that the weakest teams in the seeded groups remain stronger than the strongest teams in the unseeded groups according to the pre-tournament seeding. This design directly implies more competitive group matches, which might be important if the strengths of the teams vary substantially. Furthermore, even the lowest-ranked teams in the seeded groups qualify for the same stage in the knockout phase as the highest-ranked teams in the unseeded groups, which automatically ensures incentive compatibility.

On the other hand, in the previous design of the EHF Champions League, as well as in our proposal, some teams are eliminated from all groups. This is necessary to create appropriate incentives to perform for the top teams---otherwise, they can essentially rest in the group stage. But it also means that monotonicity should be studied carefully in order to exclude tanking by a strong team in order to be assigned to a Tier 2 group, and easily qualify for the knockout stage by playing against weak opponents in the group stage.

\subsection{The simulation model} \label{Sec32}

The outcomes of all group matches are determined by the methodology of \citet{FootballRankings2020}. The number of goals scored in a match is assumed to follow a Poisson distribution \citep{Maher1982, vanEetveldeLey2019}, and the expected number of goals is a quartic polynomial of win expectancy. The function is estimated by the least squares method based on almost 40 thousand matches between national football teams, separately for home-away games and those played on neutral ground \citep{FootballRankings2020}. Win expectancy depends on the World Football Elo Ratings (\href{http://eloratings.net/about}{http://eloratings.net/about}) and the field of the match.

Denote the expected number of goals scored by team $i$ against team $j$ by $\lambda_{ij}^{(f)}$ if the match is played on field $f$ (home: $f = h$; away: $f = a$; neutral: $f = n$). The probability that team $i$ scores $k$ goals in this match is
\begin{equation*} \label{Poisson_dist}
P_{ij}(k) = \frac{ \left( \lambda_{ij}^{(f)} \right)^k \exp \left( -\lambda_{ij}^{(f)} \right)}{k!}.
\end{equation*}

According to the World Football Elo Ratings, the win expectancy $W_{ij}$ of team $i$ with Elo $E_i$ against team $j$ with Elo $E_j$ equals
\begin{equation} \label{eq1}
W_{ij} = \frac{1}{1 + 10^{-(E_i - E_j)/400}}.
\end{equation}
Home advantage is accounted for by adding 100 points to the rating of the home team (in our case, Canada, Mexico, United States) as specified by World Football Elo Ratings (\url{https://eloratings.net/about}), which is adopted by \citet{FootballRankings2020} to calibrate the simulation model of matches played at home or away.

\citet{FootballRankings2020} report the following estimations for $\lambda_{ij}^{(f)}$ based on $W_{ij}$.
For games played on a neutral field:
\begin{equation*} \label{Exp_goals_neutral}
\lambda_{ij}^{(n)} = 
\left\{ \begin{array}{ll}
3.90388 \cdot W_{ij}^4 - 0.58486 \cdot W_{ij}^3 \\
- 2.98315 \cdot W_{ij}^2 + 3.13160 \cdot W_{ij} + 0.33193 & \textrm{if } W_{ij} \leq 0.9 \\ \\
308097.45501 \cdot (W_{ij}-0.9)^4 - 42803.04696 \cdot (W_{ij}-0.9)^3 & \\
+ 2116.35304 \cdot (W_{ij}-0.9)^2 - 9.61869 \cdot (W_{ij}-0.9) + 2.86899 & \textrm{if } W_{ij} > 0.9.
\end{array} \right. 
\end{equation*}

For games played by a host nation $H$ (Canada, Mexico, United States), the expected number of goals scored by the home team equals
\begin{equation*} \label{Exp_goals_home}
\lambda_{Hj}^{(h)} = 
\left\{ \begin{array}{ll}
-5.42301 \cdot W_{Hj}^4 + 15.49728 \cdot W_{Hj}^3 \\
- 12.6499 \cdot W_{Hj}^2 + 5.36198 \cdot W_{Hj} + 0.22863 & \textrm{if } W_{Hj} \leq 0.9 \\ \\
231098.16153 \cdot (W_{Hj}-0.9)^4 - 30953.10199 \cdot (W_{Hj}-0.9)^3 & \\
+ 1347.51495 \cdot (W_{Hj}-0.9)^2 - 1.63074 \cdot (W_{Hj}-0.9) + 2.54747 & \textrm{if } W_{Hj} > 0.9,
\end{array} \right.
\end{equation*}
while the expected number of goals scored by the away team $j$ equals
\begin{equation*} \label{Exp_goals_away}
\lambda_{Hj}^{(a)} = 
\left\{ \begin{array}{ll}
90173.57949 \cdot (W_{Hj} - 0.1)^4 + 10064.38612 \cdot (W_{Hj} - 0.1)^3 \\
+ 218.6628 \cdot (W_{Hj} - 0.1)^2 - 11.06198 \cdot (W_{Hj} - 0.1) + 2.28291 & \textrm{if } W_{Hj} < 0.1 \\ \\
-1.25010 \cdot W_{Hj}^4 -  1.99984 \cdot W_{Hj}^3 & \\
+ 6.54946 \cdot W_{Hj}^2 - 5.83979 \cdot W_{Hj} + 2.80352 & \textrm{if } W_{Hj} \geq 0.1.
\end{array} \right.
\end{equation*}
This simulation methodology has recently been applied in several studies on tournament design \citep{Csato2022a, Csato2023d, Csato2023a, Csato2023c, Csato2025c, Csato2025d, Stronka2024}.

The knockout stage consists of matches where one team qualifies and the other is eliminated. Thus, we directly use Equation~\eqref{eq1} of win expectancy. Again, the Elo rating of a host is increased by 100. Note that hosts can play against each other in the knockout stage, but the exact location of their game remains unknown; hence, the rating of both teams is increased, which is equivalent to using the original Elo ratings. The third-place game is ignored; consequently, 103 and 95 matches are simulated in the official and the imbalanced formats, respectively.

The official tie-breaking rules are implemented till the number of goals scored. Disciplinary points are ignored as yellow and red cards are not modelled, while the FIFA World Ranking is replaced by the World Football Elo Ratings, our measure of team strength.

1 million simulation runs are carried out for both the official and the imbalanced formats: one thousand random group draws are generated separately for each design, and each of them is simulated one thousand times.

The boundaries of 99\% confidence intervals for a simulated probability $p$ are given by $\pm 2.8 \sqrt{p(1-p) / n}$, and $2.8 \sqrt{p(1-p) / n} \leq 0.14\%$ due to $n = 1{,}000{,}000$. Therefore, confidence intervals will not be provided in the following because the averages of tournament metrics between the two tournament formats differ statistically significantly.

\subsection{Data} \label{Sec33}

The same set of 48 teams competes in the two formats to get comparable results. Since simulating the complicated qualification process \citep{Csato2023c} would be cumbersome, the top $\ell$ teams are selected from each confederation, similar to \citet{KrumerMoreno-Ternero2023} and \citet{Stronka2024}. There are three hosts, Canada, Mexico, and the United States. The value of $\ell$ is given by the quotas of the six confederations: AFC (8), CAF (9), CONCACAF (further 3), CONMEBOL (6), OFC (1), and UEFA (16).

Two additional teams can qualify via the inter-confederation play-offs, which will be a tournament played in one or more of the host countries. It includes six teams, one from AFC, CAF, CONMEBOL, OFC each, and two from CONCACAF. The two highest-rated seeded teams play against the winners of the first two knockout clashes between the four unseeded teams, and the two winners qualify for the 2026 FIFA World Cup group stage. In our simulations, the win expectancy Equation~\eqref{eq1} is used again to determine the winners of these matches played on neutral ground.

\begin{table}[t!]
  \centering
  \caption{National teams playing in the hypothetical 2026 FIFA World Cup}
  \label{Table1}
\centerline{
\begin{threeparttable}
    \rowcolors{1}{}{gray!20}
    \begin{tabularx}{1.1\textwidth}{lClCC c lClCC} \toprule
    Team  & Elo & Conf.\ & Pot & Pot* &       & Team  & Elo & Conf.\ & Pot   & Pot* \\ \bottomrule
    Canada & 1770  & CONC  & 1     & 1     &       & Morocco & 1780  & CAF   & 3 & 5 \\
    Mexico & 1792  & CONC  & 1     & 1     &       & Venezuela & 1771  & CONM  & 3 & 5 \\
    United States & 1727  & CONC  & 1     & 1     &       & Senegal & 1754  & CAF   & 3 & 5 \\
    Spain & 2157  & UEFA  & 1     & 1     &       & South Korea & 1743  & AFC   & 3 & 5 \\
    Argentina & 2149  & CONM  & 1     & 1     &       & Panama & 1731  & CONC  & 3 & 5 \\
    Colombia & 2061  & CONM  & 1     & 1     &       & Australia & 1723  & AFC   & 3 & 5 \\
    France & 2015  & UEFA  & 1     & 1     &       & Uzbekistan & 1690  & AFC   & 3 & 6 \\
    Brazil & 1997  & CONM  & 1     & 1     &       & Tunisia & 1681  & CAF   & 3 & 6 \\
    England & 1996  & UEFA  & 1     & 2     &       & Algeria & 1671  & CAF   & 3 & 6 \\
    Portugal & 1975  & UEFA  & 1     & 2     &       & Egypt & 1669  & CAF   & 3 & 6 \\
    Netherlands & 1960  & UEFA  & 1     & 2     &       & Costa Rica & 1662  & CONC  & 4 & 7 \\
    Germany & 1957  & UEFA  & 1     & 2     &       & Iraq  & 1658  & AFC   & 4 & 7 \\
    Uruguay & 1956  & CONM  & 2     & 2     &       & Ivory Coast & 1635  & CAF   & 4 & 7 \\
    Italy & 1943  & UEFA  & 2     & 2     &       & Jordan & 1629  & AFC   & 4 & 7 \\
    Belgium & 1923  & UEFA  & 2     & 2     &       & Jamaica & 1595  & CONC  & 4 & 7 \\
    Croatia & 1912  & UEFA  & 2     & 2     &       & Mali  & 1594  & CAF   & 4 & 7 \\
    Japan & 1873  & AFC   & 2     & 3     &       & Nigeria & 1592  & CAF   & 4 & 7 \\
    Ecuador & 1871  & CONM  & 2     & 3     &       & Saudi Arabia & 1588  & AFC   & 4 & 7 \\
    Denmark & 1863  & UEFA  & 2     & 3     &       & Angola & 1585  & CAF   & 4 & 8 \\
    Switzerland & 1855  & UEFA  & 2     & 3     &       & New Zealand & 1570  & OFC   & 4 & 8 \\
    Austria & 1853  & UEFA  & 2     & 4     &       & Peru  & 1714  & CONM  & PO    & PO \\
    Iran  & 1819  & AFC   & 2     & 4     &       & Cameroon & 1580  & CAF   & PO    & PO \\
    Turkey & 1812  & UEFA  & 2     & 4     &       & Qatar & 1574  & AFC   & PO    & PO \\
    Russia & 1785  & UEFA  & 2     & 4     &       & Haiti & 1517  & CONC  & PO    & PO \\
    Serbia & 1784  & UEFA  & 3     & 5     &       & Honduras & 1506  & CONC  & PO    & PO \\
    Greece & 1781  & UEFA  & 3     & 5     &       & New Caledonia & 1234  & OFC   & PO    & PO \\ \bottomrule
    \end{tabularx}
\begin{tablenotes} \footnotesize
\item
Abbreviations: Conf.\ = Confederation; CONC = CONCACAF; CONM = CONMEBOL.
\item
Seeding is based on the Elo ratings of 1 October 2024. Source: \url{https://www.international-football.net/elo-ratings-table?year=2024&month=10&day=01}.
\item
Russia is included only because it had the 14th highest Elo rating among UEFA teams on this day, even though it is currently suspended from FIFA competitions.
\item
The column Pot/Pot* shows the pot of the team in the official/proposed imbalanced format.
\item
In the official design, groups contain one team from Pots 1--4 each.
\item
In the imbalanced design, Tier 1 groups contain one team from Pots 1, 2, 5, 7 each, while Tier 2 groups contain one team from Pots 3, 4, 6, 8 each.
\item
The last six teams play in the inter-confederation play-offs, and the two winners are assigned to the last pot (Pot 4 or Pot 8).
\end{tablenotes}
\end{threeparttable}
}
\end{table}

\begin{figure}[t!]
\centering

\begin{tikzpicture}
\begin{axis}[
name = axis1,
xlabel = Rank of the national team,
x label style = {font=\small},
ylabel = Elo rating of the national team,
y label style = {font=\small},
width = \textwidth,
height = 0.6\textwidth,
ymajorgrids = true,
xmin = 0,
xmax = 49,
legend style = {font=\small,at={(0.275,-0.15)},anchor=north west,legend columns=2},
legend entries = {Tier 1 groups$\qquad$,Tier 2 groups}
] 
\addplot [blue, mark=asterisk, only marks, mark size=2.5pt, mark options={solid,semithick}] coordinates {
(1,1770)
(2,1792)
(3,1727)
(4,2157)
(5,2149)
(6,2061)
(7,2015)
(8,1997)
(9,1996)
(10,1975)
(11,1960)
(12,1957)
(13,1956)
(14,1943)
(15,1923)
(16,1912)
(25,1784)
(26,1781)
(27,1780)
(28,1771)
(29,1754)
(30,1743)
(31,1731)
(32,1723)
(37,1662)
(38,1658)
(39,1635)
(40,1629)
(41,1595)
(42,1594)
(43,1592)
(44,1588)
};
\addplot [ForestGreen, mark=triangle, only marks, mark size=2pt, mark options={solid,thick}] coordinates {
(17,1873)
(18,1871)
(19,1863)
(20,1855)
(21,1853)
(22,1819)
(23,1812)
(24,1785)
(33,1690)
(34,1681)
(35,1671)
(36,1669)
(45,1585)
(46,1570)
(47,1678.753)
(48,1549.413)
};
\addplot [black, mark=asterisk, only marks, mark size=2.5pt, mark options={solid,semithick}] coordinates {
(1,1870)
(2,1892)
(3,1827)
};
\end{axis}
\end{tikzpicture}

\captionsetup{justification=centerfirst}
\caption{Teams in the imbalanced format for the 2026 FIFA World Cup \\ \vspace{0.25cm}
\footnotesize{\emph{Notes}: The two winners of the play-offs are represented by their expected Elo ratings. \\
The black nodes show the three hosts if their Elo ratings are increased by 100 due to home advantage.}}
\label{Fig1}
\end{figure}


The strengths of the teams are given by their Elo ratings of 1 October 2024 and are reported in Table~\ref{Table1}.
The strengths and the allocation of the teams in the proposed imbalanced format are illustrated in Figure~\ref{Fig1}, too. Note that the 13 strongest teams and the three hosts are assigned to Tier 1 groups, that is, they can directly qualify for the Round of 16, in contrast to the next eight teams. On the other hand, the latter eight teams play against opponents that have at most 1700 Elo points in Tier 2 groups. The three hosts would be ranked 14th (Mexico), 17th (Canada), and 21st (United States) according to their strength after accounting for home advantage.

\begin{table}[t!]
  \centering
  \caption{Matches in the last round of the group stage in the recent FIFA World Cups}
  \label{Table2}
\begin{threeparttable}
    \rowcolors{1}{}{gray!20}
    \begin{tabularx}{0.75\textwidth}{lCC} \toprule
    Group & 2018 FIFA World Cup & 2022 FIFA World Cup \\ \bottomrule
    A     & 1v2, 3v4 & 1v2, 3v4 \\
    B     & 1v3, 2v4 & 1v4, 2v3 \\
    C     & 1v3, 2v4 & 1v3, 2v4 \\
    D     & 1v4, 2v3 & 1v3, 2v4 \\
    E     & 1v3, 2v4 & 1v3, 2v4 \\
    F     & 1v4, 2v3 & 1v2, 3v4 \\
    G     & 1v2, 3v4 & 1v4, 2v3 \\
    H     & 1v4, 2v3 & 1v3, 2v4 \\ \bottomrule
    \end{tabularx}
\begin{tablenotes} \footnotesize
\item
\emph{Note}: $k$v$\ell$ means that the team drawn from Pot $k$ plays against the team drawn from Pot $\ell$ in the last round, $1 \leq k, \ell \leq 4$.
\end{tablenotes}
\end{threeparttable}
\end{table}

It is well-known that the probability of stakeless matches is influenced by the schedule of the matches \citep{ChaterArrondelGayantLaslier2021, CsatoMolontayPinter2024, Guyon2020a, Stronka2024}, even if our simulation method is independent of the order in which the games are played. In both the 2018 \citep{FIFA2017c} and 2022 FIFA World Cups \citep{FIFA2022a}, a random schedule was used since, after the group of a team was determined, the position in which it plays depended on an additional random draw.
As will be discussed in Section~\ref{Sec34}, it is sufficient to focus on the set of matches played simultaneously in the last round, where three options exist, distinguished by the opponent of the team drawn from Pot 1. Table~\ref{Table2} shows the actual schedules of the eight groups in the two recent FIFA World Cups; each of the three possible schedules occurred in at least two groups in both editions.

\subsection{Evaluation metrics} \label{Sec34}

The two alternative formats are assessed primarily on the basis of the likelihood that a group match provides no incentives for the teams to perform. Previous studies analyse stakeless matches, where a team is indifferent to the outcome \citep{ChaterArrondelGayantLaslier2021, CsatoMolontayPinter2024, Gyimesi2024}. \citet{CsatoMolontayPinter2024} distinguish weakly (if one team is indifferent) and strongly (if both teams are indifferent) stakeless matches. However, none of these studies deals with the strength or the current position of the opposing teams. To take these aspects into account, stakeless matches are considered at the level of teams, which is our main methodological innovation.

Denote the number of stakeless group matches of team $i$ by $S_i$.
$S_i$ is increased by one if team $i$ plays a match whose outcome does not affect the probability of any of the following events:
\begin{itemize}
\item 
Being eliminated;
\item 
Qualification for the Round of 32 (in the official format) or the knockout stage play-offs (in the proposed imbalanced format);
\item
Qualification for the Round of 16 (in the proposed imbalanced format).
\end{itemize}
Obviously, only the last round can contain a stakeless match if four teams play a single round-robin competition, since a team could have at most three points before its second match, and six points can be collected in two matches. In addition, the last round is played simultaneously since the Disgrace of Gij\'on \citep[Section~3.1]{KendallLenten2017}. Thus, the average value of $S_i$ over all simulation runs gives the probability that the match played by team $i$ in the last round of the group stage becomes stakeless.

The status of the team that plays a stakeless match could also be important. An already eliminated team might still exert full effort in a stakeless match as they have no reason to rest their players. On the other hand, an already qualified team has a powerful incentive to rest and reserve energy. Therefore, stakeless games can be classified into two categories: matches played by teams already qualified ($S_i^A$) and matches played by teams already eliminated ($S_i^E$). 

$S_i^A$ is higher if more teams advance from a group to the knockout stage. Consequently, it makes a substantial difference in the official format whether two or three teams qualify from a certain group. We provide a minimum and a maximum value of $S_i^A$ accordingly. In the imbalanced format, the minimum and the maximum coincide, as the number of teams advancing is deterministic and does not depend on the results in other groups.

In the official format, the maximum $S_i^{A,\max}$ is reached if the match played by a team that is guaranteed to be in the top three in their group before the last round is called stakeless. However, it does not happen in practice in all groups, only in those where some teams already know that they would certainly advance even if they are ranked third.
Contrarily, the minimum $S_i^{A,\min}$ is reached if no team has any information on whether the third position would be sufficient for them to qualify. This also happens in some, but usually not all, groups. Therefore, the true value $S_i^A$ lies somewhere between $S_i^{A,\min}$ and $S_i^{A,\max}$.

The same logic applies to the teams that have already been eliminated. In the imbalanced format, there is no difference between $S_i^{E,\min}$ and $S_i^{E,\max}$. In the official format, $S_i^{E,\max}$ applies if the third place is insufficient to qualify for the Round of 32, and $S_i^{E,\min}$ applies if it might mean qualification.

As we have argued in the Introduction, a stakeless match is less costly if the team without any incentive to perform has \emph{a priori} a lower chance of winning the match. A reasonable weighting factor is the win expectancy $W_{ij}$ according to Equation~\eqref{eq1}. Hence, the expected number of weighted stakeless matches is calculated as
\[
S_i^W =  S_i \cdot W_{ij}.
\]
For instance, the stakeless matches discussed in Section~\ref{Sec1} are associated with the weights 0.8835 (the winning probability of France against Tunisia because France is already qualified) and 0.3351 (the winning probability of Canada against Morocco because Canada is already eliminated), respectively.

However, $S_i$ is not properly defined in the official format because the number of teams qualified for the Round of 32 can be either two or three in a group. Therefore, the \emph{minimum} of the expected number of weighted stakeless matches is computed as a lower bound by assuming that the teams do not know in advance whether the third place will be sufficient for qualification, namely, three different prizes exist (top 2, third place, fourth place), and a match is stakeless for a team only if its outcome has no effect on the probability of obtaining these three prizes. This is the best case with respect to the competitiveness of the matches since the third-ranked team is guaranteed to know its target to qualify in the last group because the results of the previous groups become known. Therefore, similar to \citet[Section~5.2]{ChaterArrondelGayantLaslier2021}, we also compute weighted stakeless matches in this worst case scenario.

\subsection{Limitations} \label{Sec35}

Naturally, a careful interpretation of our results is warranted.
We consider only one set of Elo ratings and a different distribution of strengths may change the quantitative findings. This is less sophisticated than the approach of \citet{Stronka2024} who takes these values from three different dates. 
Second, the match outcomes are assumed to be independent of the schedule, even though the economic literature has found robust evidence for the effect of the order of the games \citep{LaicaLauberSahm2021}. In particular, playing in the first and third matches of the group stage leads to a significantly higher probability of qualification in the FIFA World Cup \citep{KrumerLechner2017}.
Third, only one reasonable variant of the imbalanced format is studied, however, the composition of Tier 1 and Tier 2 groups is not straightforward.
Fourth, World Football Elo Ratings have been adopted as the measure of team strength, albeit the official FIFA World Ranking also uses the Elo system since 2018 \citep{FIFA2018c} and could help persuade the decision-makers. On the other hand, the current FIFA World Ranking still suffers from several shortcomings, such as the counterproductive weighting of game importance, the absence of home advantage, and ignoring goal difference \citep{SzczecinskiRoatis2022}.
Fifth, the home advantage enjoyed by the host nations is accounted for in a relatively simple framework; however, the proportion of group matches played by them (9 out of 72, 12.5\%) remains low, hence, the home advantage parameter does not have a large effect on our results.
Sixth, regional rivalries might provide extra incentives even if a match seems to be stakeless in the proposed model. On the other hand, a group can contain at most one team from one confederation except for UEFA (if the two, rather weak winners of the play-offs are not considered), and the composition of the pots substantially decreases the relevance of regional rivalries even for UEFA teams; for example, France cannot play against Spain, and England cannot play against Germany due to the seeding shown in Table~\ref{Table1}.
Finally, the regression model of \citet{FootballRankings2020} may be suboptimal for the FIFA World Cup where the weakest national teams have never played. 

All probabilities are derived using Monte Carlo simulations.
However, the schedule of the FIFA World Cup can be exploited to compute the winning probabilities exactly: \citet{BrandesMarmullaSmokovic2025} give an algorithm that is two orders of magnitude faster than any reasonably accurate approximation. It remains to be seen whether this result could be used in our case with 48 teams, one thousand different group draws, and knockout brackets randomised subject to the basic constraints.

To conclude, even though our Elo-based approach is not necessarily the best available simulation model, it is mainly used to choose between two competition formats. Thus, the counterargument of \citet[p.~534]{Appleton1995} applies: ``\emph{Since our intention is to compare tournament designs, and not to estimate the chance of the best player winning any particular tournament, we may within reason take whatever model of determining winners that we please}''.

\section{Results} \label{Sec4}

Section~\ref{Sec41} evaluates the aggregated measures of stakeless games both in the official and the proposed tournament formats for the 2026 FIFA World Cup, as well as discusses monotonicity, which is a crucial issue in the presence of imbalanced groups. After that, Section~\ref{Sec42} compares the two designs from several aspects at the level of individual teams.

\subsection{Aggregated measures of stakeless matches and monotonicity} \label{Sec41}

\begin{table}[t!]
  \centering
  \caption{Aggregated tournament metrics for the two formats}
  \label{Table3}
\begin{threeparttable}
    \rowcolors{1}{}{gray!20}
    \begin{tabularx}{0.8\textwidth}{l CC} \toprule
    Format & Official & Imbalanced \\ \bottomrule
    Average Elo difference & 214.4 & 208.2 \\ 
    Average $S_i^{A,\min}$ & 0.198 & 0.072 \\
    Average $S_i^{A,\max}$ & 0.434 & 0.072 \\
    Average $S_i^{E,\min}$ & 0.027 & 0.159 \\
    Average $S_i^{E,\max}$ & 0.166 & 0.159 \\
    $S_i^W(\min)$ & 0.129 & 0.121 \\
    $S_i^W(22)$ & 0.193 & --- \\
    $S_i^W(31)$ & 0.268 & --- \\ \toprule
    \end{tabularx}
\begin{tablenotes} \footnotesize
\item
$S_i^{A,\min}$ ($S_i^{A,\max}$) assumes that two (three) teams qualify in the official format.
\item
$S_i^{E,\min}$ ($S_i^{E,\max}$) assumes that three (two) teams qualify in the official format.
\item
$S_i^W(\min)$ conservatively assumes that three prizes exist in the official format, see Section~\ref{Sec34}.
\item 
$S_i^W(22)$ ($S_i^W(31)$) assumes that two (three) teams qualify in the official format.
\end{tablenotes}
\end{threeparttable}
\end{table}

According to Table~\ref{Table3}, the outcome of group stage matches is, on average, more uncertain in the modified format than in the original, as the difference between the Elo ratings is smaller. Furthermore, the imbalanced format provides a substantially lower ratio of stakeless matches for teams that are qualified after two rounds, even if the best case is assumed in the official format.
On the other hand, the probability of stakeless matches played by eliminated teams does not decrease compared to the official format. This is not surprising since not only 16 but 24 teams are eliminated at the end of the group stage.

Last but not least, the imbalanced format is clearly better than even the best case in the official format if the stakeless matches are weighted according to the winning probabilities as suggested in Section~\ref{Sec34}. Its advantage becomes much stronger for the last group(s), where the teams already know whether the third place would be sufficient for qualification or not. Thus, the aggregated measures of stakeless matches tend to favour the proposed imbalanced design.

Despite these benefits, the imbalanced design may create misaligned incentives if a team could be better off by being assigned to a weaker pot \citep{Csato2020b, Csato2021a}, which is impossible in the official design. This is far from clear as, for example, the teams drawn from Pot 4 play against weaker opponents in Tier 2 groups where direct qualification for the Round of 16 is excluded, but the teams drawn from Pot 5 play against stronger opponents in Tier 1 groups where they can directly qualify for the Round of 16.

\begin{figure}[t!]
\centering

\begin{tikzpicture}
\begin{axis}[
name = axis1,
xlabel = Elo rating of the national team,
x label style = {font=\small},
ylabel = Probability of qualifying for the Round of 16,
y label style = {font=\small},
width = \textwidth,
height = 0.6\textwidth,
ymajorgrids = true,
ymin = 0,
ymax = 1.02,
max space between ticks = 50,
legend style = {font=\small,at={(0,-0.15)},anchor=north west,legend columns=3},
legend entries = {Official format$\qquad$, {Imbalanced, Tier 1 groups}$\qquad$,{Imbalanced, Tier 2 groups}}
]
\addplot [red, mark=oplus, only marks, mark size=2pt, mark options={solid,thin}] coordinates {
(1792,0.518601)
(1770,0.479145)
(1727,0.400777)
(2157,0.862801)
(2149,0.857258)
(2061,0.775415)
(2015,0.714999)
(1997,0.688319)
(1996,0.685734)
(1975,0.654028)
(1960,0.63149)
(1957,0.625807)
(1956,0.600216)
(1943,0.582801)
(1923,0.543804)
(1912,0.524888)
(1873,0.443501)
(1871,0.444681)
(1863,0.431732)
(1855,0.41548)
(1853,0.414703)
(1819,0.341674)
(1812,0.336767)
(1785,0.29027)
(1784,0.255343)
(1781,0.253319)
(1780,0.241196)
(1771,0.23878)
(1754,0.201173)
(1743,0.19125)
(1731,0.166462)
(1723,0.162674)
(1690,0.124235)
(1681,0.112879)
(1671,0.101252)
(1669,0.099919)
(1662,0.072776)
(1658,0.074546)
(1635,0.056515)
(1629,0.055237)
(1595,0.034911)
(1594,0.035459)
(1592,0.0354)
(1588,0.034925)
(1585,0.032419)
(1570,0.028868)
(1678.753,0.101123)
(1549.413,0.024448)
}; 
\addplot [blue, mark=asterisk, only marks, mark size=2.5pt, mark options={solid,semithick}] coordinates {
(1792,0.604386)
(1770,0.562603)
(1727,0.476946)
(2157,0.953438)
(2149,0.948647)
(2061,0.873168)
(2015,0.811056)
(1997,0.781954)
(1996,0.752029)
(1975,0.716183)
(1960,0.688712)
(1957,0.682294)
(1956,0.668314)
(1943,0.656534)
(1923,0.618307)
(1912,0.597535)
(1784,0.216265)
(1781,0.208416)
(1780,0.183018)
(1771,0.210947)
(1754,0.149643)
(1743,0.136769)
(1731,0.095539)
(1723,0.112578)
(1662,0.037613)
(1658,0.053697)
(1635,0.044383)
(1629,0.039201)
(1595,0.015725)
(1594,0.026181)
(1592,0.026592)
(1588,0.023885)
};
\addplot [ForestGreen, mark=triangle, only marks, mark size=2pt, mark options={solid,thick}] coordinates {
(1873,0.378518)
(1871,0.374212)
(1863,0.360768)
(1855,0.347422)
(1853,0.33521)
(1819,0.282178)
(1812,0.26916)
(1785,0.230747)
(1690,0.088092)
(1681,0.078815)
(1671,0.071297)
(1669,0.06963)
(1585,0.023914)
(1570,0.020367)
(1678.753,0.079794)
(1549.413,0.017318)
};
\addplot [black, mark=oplus, only marks, mark size=2pt, mark options={solid,thin}] coordinates {
(1892,0.518601)
(1870,0.479145)
(1827,0.400777)
};
\addplot [black, mark=asterisk, only marks, mark size=2.5pt, mark options={solid,semithick}] coordinates {
(1892,0.604386)
(1870,0.562603)
(1827,0.476946)
};
\end{axis}
\end{tikzpicture}

\captionsetup{justification=centerfirst}
\caption{The probability of qualification for the Round of 16 in the two formats for the 2026 FIFA World Cup \\ \vspace{0.25cm}
\footnotesize{\emph{Notes}: The two winners of the play-offs are represented by their expected Elo ratings. \\
The black nodes show the three hosts if their Elo ratings are increased by 100 due to home advantage.}}
\label{Fig2}
\end{figure}


Therefore, Figure~\ref{Fig2} plots the chance of advancing to the Round of 16 as a function of the Elo rating. No anomaly can be seen, that is, teams close to the boundary of two pots cannot benefit from being in a lower-ranked pot. The three host nations are favoured by their guaranteed place in Pot 1 as the black nodes clearly lie above the trend line---but this is a deliberate policy of the organiser. 
The imbalanced design is better for the three hosts and the 13 strongest teams, while it punishes the other teams compared to the official format.

Three national teams (Panama, Costa Rica, Jamaica) have a lower chance to qualify for the Round of 16 than teams of similar Elo ratings from other confederations, especially in the imbalanced format. The reason is clear: they cannot play in the groups of the hosts that are weaker than the other teams drawn from Pot 1. This effect is stronger in the suggested design, which contains eight Tier 1 groups instead of the 12 in the official design. In contrast, Venezuela, the weakest CONMEBOL team, benefits from not playing against the strong South American teams.
This observation calls our attention to the non-negligible sporting effects of draw constraints \citep{Csato2022d} that are stronger with fewer groups.

\begin{table}[t!]
  \centering
  \caption{The maximal possible gain from tanking in the imbalanced design}
  \label{Table4}
\begin{threeparttable}
    \rowcolors{1}{gray!20}{}
    \begin{tabularx}{\textwidth}{l cc RRR} \toprule \hiderowcolors
    \multirow{2}[0]{*}{Team} & \multicolumn{2}{c}{Pot} & \multicolumn{3}{c}{Probability of qualifying for the Round of 16} \\
          & Original & Changed & Original & Changed & Difference \\ \bottomrule \showrowcolors
    Spain & 1     & 3     & 95.34\% & 80.29\% & $-$15.05\% \\
    England & 2     & 3     & 75.20\% & 60.04\% & $-$15.16\% \\
    Japan & 3     & 5     & 37.85\% & 38.37\% & 0.52\% \\
    Austria & 4     & 5     & 33.52\% & 36.25\% & 2.73\% \\
    Serbia & 5     & 6     & 21.63\% & 18.93\% & $-$2.70\% \\
    Uzbekistan & 6     & 7     & 8.81\% & 5.06\% & $-$3.75\% \\
    Costa Rica & 7     & 8     & 3.76\% & 5.73\% & 1.97\% \\ \toprule
    \end{tabularx}
\begin{tablenotes} \footnotesize
\item
The column Changed pot indicates the strongest pot from which the groups of the other tier are drawn that is available by tanking.
\end{tablenotes}
\end{threeparttable}
\end{table}

But the monotonicity of Figure~\ref{Fig2} does not imply the absence of misaligned incentives. In order to investigate this issue in more detail, we have considered whether the strongest teams of Pots 1--7 would be better off if they became the highest-ranked team of the strongest pot corresponding to the other tier of groups due to a successful manipulation strategy, while all other teams are assigned according to their true strength.
The results are reported in Table~\ref{Table4}. Teams in the three strongest pots (1, 2, 5) associated with Tier 1 groups could not benefit from losing, similar to the two weakest pots (6, 8) associated with Tier 2 groups. Consequently, the 16 strongest teams have the right incentives in our proposal.

On the other hand, the strongest teams of Pots 3, 4, and 7 could profit from a successful tanking. However, the potential gain always remains below three percentage points in terms of the probability of reaching the Round of 16, which makes the consideration of such a risky and uncertain strategy improbable. The payoff of a manipulation might be reduced by modifying the allocation of pots to Tiers 1 and 2. Another solution could be adding more uncertainty to the group draw, for example, by randomising the assignment of the teams to the pots, analogous to the draft lottery of the NBA and the NHL \citep{GerchakMausserMagazine1995, Schmidt2024}.

To conclude, these results indicate that the analysis of monotonicity is a crucial issue in a tournament format containing imbalanced groups, and should be the topic of further research due to the robust advantages of this competition design.

\subsection{Detailed comparison of the two tournament formats} \label{Sec42}

\begin{figure}[t!]
\centering

\begin{tikzpicture}
\begin{axis}[
name = axis1,
xlabel = Elo rating of the national team,
x label style = {font=\small},
ylabel = Average number of matches played,
y label style = {font=\small},
width = \textwidth,
height = 0.6\textwidth,
ymajorgrids = true,
max space between ticks = 50,
legend style = {font=\small,at={(0,-0.15)},anchor=north west,legend columns=3},
legend entries = {Official format$\qquad$, {Imbalanced, Tier 1 groups}$\qquad$,{Imbalanced, Tier 2 groups}}
]
\addplot [red, mark=oplus, only marks, mark size=2pt, mark options={solid,thin}] coordinates {
(1792,4.809625)
(1770,4.696925)
(1727,4.484892)
(2157,6.442865)
(2149,6.389162)
(2061,5.835751)
(2015,5.53299)
(1997,5.423839)
(1996,5.400636)
(1975,5.27576)
(1960,5.192108)
(1957,5.168045)
(1956,5.109893)
(1943,5.045677)
(1923,4.913754)
(1912,4.854835)
(1873,4.610735)
(1871,4.619884)
(1863,4.574734)
(1855,4.533912)
(1853,4.525608)
(1819,4.322326)
(1812,4.315757)
(1785,4.183496)
(1784,4.064182)
(1781,4.057984)
(1780,4.015672)
(1771,4.018265)
(1754,3.901136)
(1743,3.873839)
(1731,3.789244)
(1723,3.786149)
(1690,3.660073)
(1681,3.614464)
(1671,3.576059)
(1669,3.569208)
(1662,3.421368)
(1658,3.432605)
(1635,3.355471)
(1629,3.355243)
(1595,3.258647)
(1594,3.260738)
(1592,3.2609)
(1588,3.263904)
(1585,3.248283)
(1570,3.232609)
(1678.753,3.518313)
(1549.413,3.202435)
}; 
\addplot [blue, mark=asterisk, only marks, mark size=2.5pt, mark options={solid,semithick}] coordinates {
(1792,4.384222)
(1770,4.278407)
(1727,4.080314)
(2157,5.871783)
(2149,5.826921)
(2061,5.345427)
(2015,5.060376)
(1997,4.947468)
(1996,4.958132)
(1975,4.834265)
(1960,4.739792)
(1957,4.731269)
(1956,4.714433)
(1943,4.639666)
(1923,4.531734)
(1912,4.467329)
(1784,3.540203)
(1781,3.524359)
(1780,3.472113)
(1771,3.522789)
(1754,3.394858)
(1743,3.366788)
(1731,3.279073)
(1723,3.310545)
(1662,3.120617)
(1658,3.154865)
(1635,3.133177)
(1629,3.118637)
(1595,3.059335)
(1594,3.086239)
(1592,3.087209)
(1588,3.078461)
};
\addplot [ForestGreen, mark=triangle, only marks, mark size=2pt, mark options={solid,thick}] coordinates {
(1873,4.394384)
(1871,4.382429)
(1863,4.342343)
(1855,4.303115)
(1853,4.262634)
(1819,4.110218)
(1812,4.06971)
(1785,3.959547)
(1690,3.456133)
(1681,3.417665)
(1671,3.390135)
(1669,3.383435)
(1585,3.17358)
(1570,3.156262)
(1678.753,3.400692)
(1549.413,3.136912)
};
\addplot [black, mark=oplus, only marks, mark size=2pt, mark options={solid,thin}] coordinates {
(1892,4.809625)
(1870,4.696925)
(1827,4.484892)
};
\addplot [black, mark=asterisk, only marks, mark size=2.5pt, mark options={solid,semithick}] coordinates {
(1892,4.384222)
(1870,4.278407)
(1827,4.080314)
};
\end{axis}
\end{tikzpicture}

\captionsetup{justification=centerfirst}
\caption{The expected number of matches in the two formats for the 2026 FIFA World Cup \\ \vspace{0.25cm}
\footnotesize{\emph{Notes}: The two winners of the play-offs are represented by their expected Elo ratings. \\
The black nodes show the three hosts if their Elo ratings are increased by 100 due to home advantage.}}
\label{Fig3}
\end{figure}


The proposed imbalanced format contains the same number of group matches as the official design ($12 \times 6 = 72$) but only 24 knockout games instead of 32. The average number of matches played by each team is presented in Figure~\ref{Fig3}. All teams play fewer matches, especially the strongest teams and the eight teams assigned to Tier 1 groups from Pot 5, which are more likely to be eliminated in the group stage.

The decreased workload is favourable since top level players should play many matches, and they face a serious risk of sustaining injuries during matches \citep{Lopez-ValencianoRuiz-PerezGarcia-GomezVera-GarciaCroixMyerAyala2020}, especially during congested periods of matches \citep{BengtssonEkstrandHagglund2013, DellalLago-PenasReyChamariOrhant2015, EkstrandWaldenHagglund2004}. Tactical performance may also worsen in congested fixtures as the movement of the players becomes more time synchronised \citep{FolgadoDuarteMarquesSampaio2015}.

\begin{figure}[t!]
\centering

\begin{tikzpicture}
\begin{axis}[
name = axis1,
xlabel = Elo rating of the national team,
x label style = {font=\small},
ylabel = Average difference in Elo ratings,
y label style = {font=\small},
width = \textwidth,
height = 0.6\textwidth,
ymajorgrids = true,
max space between ticks = 50,
legend style = {font=\small,at={(0,-0.15)},anchor=north west,legend columns=3},
legend entries = {Official format$\qquad$, {Imbalanced, Tier 1 groups}$\qquad$,{Imbalanced, Tier 2 groups}}
]
\addplot [red, mark=oplus, only marks, mark size=2pt, mark options={solid,thin}] coordinates {
(1792,161.941666666667)
(1770,144.754666666667)
(1727,120.905666666667)
(2157,417.424)
(2149,413.059)
(2061,323.629)
(2015,276.455333333333)
(1997,261.105333333333)
(1996,255.856333333333)
(1975,236.411666666667)
(1960,222.559333333333)
(1957,217.990666666667)
(1956,208.973)
(1943,213.493666666667)
(1923,201.999)
(1912,195.051666666667)
(1873,172.887333333333)
(1871,174.639666666667)
(1863,171.906333333333)
(1855,169.797666666667)
(1853,166.311333333333)
(1819,153.644666666667)
(1812,153.598)
(1785,142.457)
(1784,150.867)
(1781,148.506666666667)
(1780,155.753)
(1771,148.489)
(1754,164.783)
(1743,170.142666666667)
(1731,187.318)
(1723,179.542333333333)
(1690,189.313)
(1681,190.52)
(1671,195.471)
(1669,195.547666666667)
(1662,211.754666666667)
(1658,202.468666666667)
(1635,230.634666666667)
(1629,231.924333333333)
(1595,279.487)
(1594,271.938333333333)
(1592,271.875666666667)
(1588,272.233666666667)
(1585,277.426333333333)
(1570,288.434333333333)
(1678.753,187.610333333333)
(1549.413,311.688)
}; 
\addplot [blue, mark=asterisk, only marks, mark size=2.5pt, mark options={solid,semithick}] coordinates {
(1792,154.243666666667)
(1770,147.143)
(1727,132.118)
(2157,381.223666666667)
(2149,372.960666666667)
(2061,285.134333333333)
(2015,239.102333333333)
(1997,221.507)
(1996,238.293333333333)
(1975,224.643666666667)
(1960,217.167)
(1957,215.964)
(1956,216.276666666667)
(1943,206.721666666667)
(1923,194.962333333333)
(1912,188.737333333333)
(1784,169.953666666667)
(1781,172.648666666667)
(1780,186.612333333333)
(1771,163.366666666667)
(1754,192.81)
(1743,200.847333333333)
(1731,228.105)
(1723,208.989666666667)
(1662,267.749)
(1658,239.132666666667)
(1635,254.496333333333)
(1629,267.026666666667)
(1595,333.461333333333)
(1594,298.720333333333)
(1592,297.558333333333)
(1588,307.855333333333)
};
\addplot [ForestGreen, mark=triangle, only marks, mark size=2pt, mark options={solid,thick}] coordinates {
(1873,176.665)
(1871,173.562666666667)
(1863,165.857666666667)
(1855,158.192666666667)
(1853,148.885)
(1819,135.188333333333)
(1812,134.463)
(1785,127.075)
(1690,135.473666666667)
(1681,139.747333333333)
(1671,145.389666666667)
(1669,145.201333333333)
(1585,205.473666666667)
(1570,215.172333333333)
(1678.753,127.222333333333)
(1549.413,236.499666666667)
};
\addplot [black, mark=oplus, only marks, mark size=2pt, mark options={solid,thin}] coordinates {
(1892,161.941666666667)
(1870,144.754666666667)
(1827,120.905666666667)
};
\addplot [black, mark=asterisk, only marks, mark size=2.5pt, mark options={solid,semithick}] coordinates {
(1892,154.243666666667)
(1870,147.143)
(1827,132.118)
};
\end{axis}
\end{tikzpicture}

\captionsetup{justification=centerfirst}
\caption{The uncertainty of group stage matches in the two formats for the 2026 FIFA World Cup for each national team \\ \vspace{0.25cm}
\footnotesize{\emph{Notes}: The two winners of the play-offs are represented by their expected Elo ratings. \\
The six black nodes show the three hosts if their Elo ratings are increased by 100 due to home advantage.}}
\label{Fig4}
\end{figure}


According to Figure~\ref{Fig4}, the strongest teams and 12 weak teams assigned to Tier 2 groups---from Pots 4, 6, 8---play more interesting matches in the imbalanced format. This is important because the top teams can win their first two matches more easily (such as Brazil, France, or Portugal in the 2022 FIFA World Cup) or qualify even by losing their first match without any consequence (such as Argentina in the 2022 FIFA World Cup).
Furthermore, emerging football nations have a lower chance of suffering humiliating defeats.

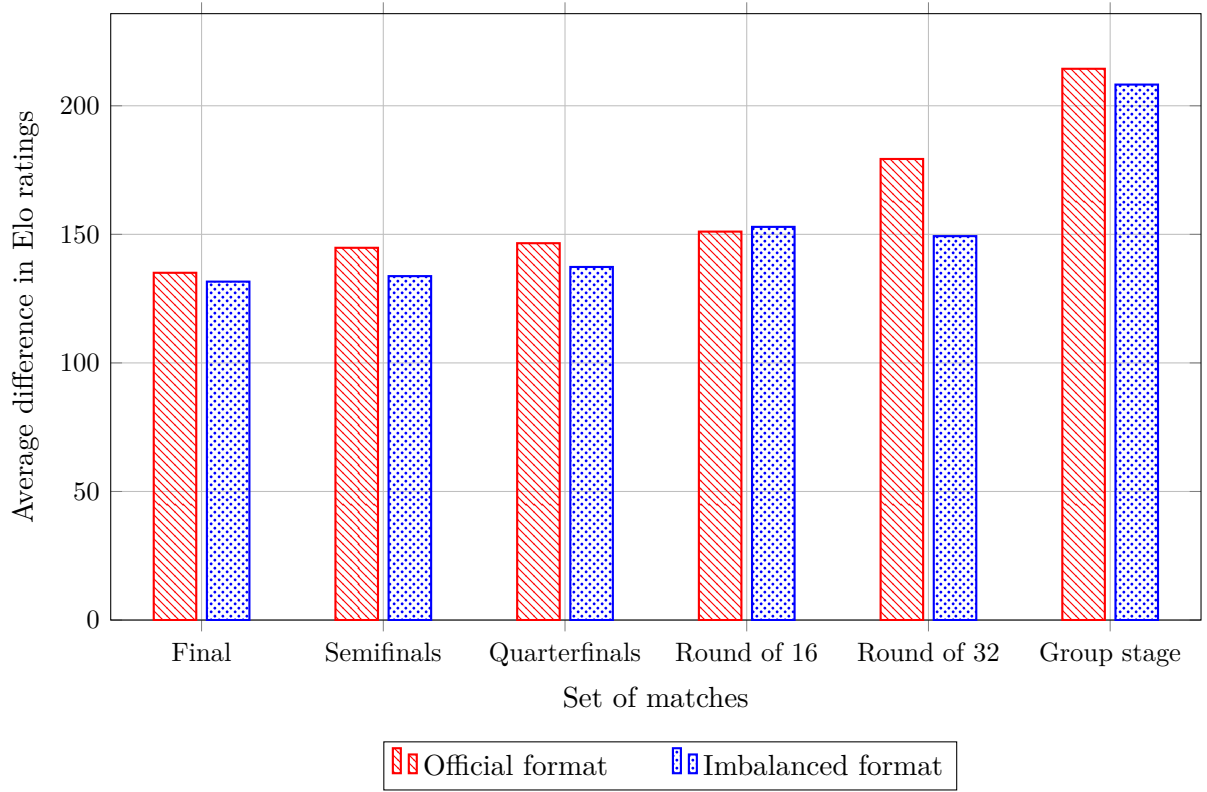
\begin{figure}[t!]
\centering

\begin{tikzpicture}
\begin{axis}[
width = \textwidth, 
height = 0.6\textwidth,
xmajorgrids,
ymajorgrids,
scaled x ticks = false,
xlabel = {Set of matches},
xlabel style = {align=center, font=\small},
xtick = data,
symbolic x coords = {Final,Semifinals,Quarterfinals,Round of 16,Round of 32,Group stage},
enlarge x limits = 0.1,
ylabel = {Average difference in Elo ratings},
ylabel style = {align=center, font=\small},
ymin = 0,
ybar = 4pt,
bar width = 16pt,
legend style = {font=\small,at={(0.25,-0.2)},anchor=north west,legend columns=3},
legend entries = {Official format$\qquad$, Imbalanced format},
]
\addplot [red, thick, pattern = north west lines, pattern color = red] coordinates{
(Final,135.049046)
(Semifinals,144.768557)
(Quarterfinals,146.5538575)
(Round of 16,151.065837875)
(Round of 32,179.310435875)
(Group stage,214.387111111111)
};
\addplot [blue, thick, pattern = crosshatch dots, pattern color = blue] coordinates{
(Final,131.612014)
(Semifinals,133.7483365)
(Quarterfinals,137.3121203)
(Round of 16,152.905014)
(Round of 32,149.3049646)
(Group stage,208.2416944)
};
\end{axis}
\end{tikzpicture}

\captionsetup{justification=centerfirst}
\caption{The uncertainty of matches in the two formats for the 2026 FIFA World Cup \\ \vspace{0.25cm}
\footnotesize{\emph{Note}: Round of 32 in the imbalanced format refers to the play-offs for the Round of 16.}}
\label{Fig5}

\end{figure}


Figure~\ref{Fig5} extends the analysis of match uncertainty to the knockout stage. The imbalanced format implies closer matches in all rounds except for the Round of 16. This is probably caused by the more efficient pairing of stronger versus weaker teams in the Round of 16, where the group winners of Tier 1 groups cannot play against each other. Hence, substantial benefits can be seen in the following two rounds, in the semifinals and the quarterfinals. The measures for the Round of 32 cannot be compared because the imbalanced format contains only eight games here.
To summarise, our proposal results in more exciting matches, even in the knockout phase.

\begin{figure}[t!]
\centering

\begin{tikzpicture}
\begin{axis}[
name = axis1,
xlabel = Value of $k$,
x label style = {font=\small},
ylabel = Ratio of matches between the $k$ strongest teams,
y label style = {font=\small},
width = \textwidth,
height = 0.6\textwidth,
xmin = 1,
xmax = 25,
ymajorgrids = true,
ymin = 0,
max space between ticks = 50,
legend style = {font=\small,at={(0.05,-0.15)},anchor=north west,legend columns=2},
legend entries = {{Official format, all matches}$\quad \, \, \qquad$, {Imbalanced format, all matches}$\quad \, \,$, {Official format, group matches}$\qquad$, {Imbalanced format, group matches}},
]
\addplot [red, mark=oplus, only marks, mark size=2pt, mark options={solid,thin}] coordinates {
(2,0.0027692427184466)
(3,0.00665174757281553)
(4,0.0110341747572816)
(5,0.0160933106796117)
(6,0.0221893883495146)
(7,0.0286830097087379)
(8,0.035526067961165)
(9,0.0430523203883495)
(10,0.0597567766990291)
(11,0.0749211650485437)
(12,0.0901267184466019)
(13,0.105509796116505)
(14,0.117446990291262)
(15,0.134772834951456)
(16,0.152710815533981)
(17,0.166673970873786)
(18,0.184926029126214)
(19,0.203410650485437)
(20,0.222161669902913)
(21,0.239311077669903)
(22,0.258443786407767)
(23,0.277796932038835)
(24,0.296449058252427)
}; 
\addplot [blue, mark=asterisk, only marks, mark size=2.5pt, mark options={solid,semithick}] coordinates {
(2,0.00336047368421053)
(3,0.00798355789473684)
(4,0.0131145157894737)
(5,0.0188637473684211)
(6,0.0313142315789474)
(7,0.0441895157894737)
(8,0.0571126842105263)
(9,0.0709238421052632)
(10,0.0853276)
(11,0.0995741157894737)
(12,0.113441231578947)
(13,0.1275662)
(14,0.145464768421053)
(15,0.156048357894737)
(16,0.166565473684211)
(17,0.185244473684211)
(18,0.196064357894737)
(19,0.206568736842105)
(20,0.227290873684211)
(21,0.246192705263158)
(22,0.266067642105263)
(23,0.285592821052632)
(24,0.304147189473684)
};
\addplot [orange, mark=otimes, only marks, mark size=2pt, mark options={solid,thin}] coordinates {
(2,0)
(3,0)
(4,0)
(5,0)
(6,0)
(7,0)
(8,0)
(9,0)
(10,0.0122777777777778)
(11,0.0219166666666667)
(12,0.0314861111111111)
(13,0.0409305555555556)
(14,0.0454444444444444)
(15,0.0581527777777778)
(16,0.071)
(17,0.0772638888888889)
(18,0.0896388888888889)
(19,0.102083333333333)
(20,0.114541666666667)
(21,0.125)
(22,0.138888888888889)
(23,0.152777777777778)
(24,0.166666666666667)
}; 
\addplot [ForestGreen, mark=star, only marks, mark size=2.5pt, mark options={solid,semithick}] coordinates {
(2,0)
(3,0)
(4,0)
(5,0)
(6,0.00851388888888889)
(7,0.0172777777777778)
(8,0.0256805555555556)
(9,0.034625)
(10,0.0437361111111111)
(11,0.0522638888888889)
(12,0.0607083333333333)
(13,0.0694444444444444)
(14,0.0833333333333333)
(15,0.0833333333333333)
(16,0.0833333333333333)
(17,0.0972222222222222)
(18,0.0972222222222222)
(19,0.0972222222222222)
(20,0.111111111111111)
(21,0.125)
(22,0.138888888888889)
(23,0.152777777777778)
(24,0.166666666666667)
};
\end{axis}
\end{tikzpicture}

\captionsetup{justification=centerfirst}
\caption{The ratio of matches between the $k$ strongest teams in the two formats for the 2026 FIFA World Cup \\ \vspace{0.25cm}
\footnotesize{\emph{Notes}: The three hosts are counted according to their Elo rating increased by 100 due to home advantage. \\
The two winners of the play-offs are counted as the 47th and 48th teams.}}
\label{Fig6}
\end{figure}
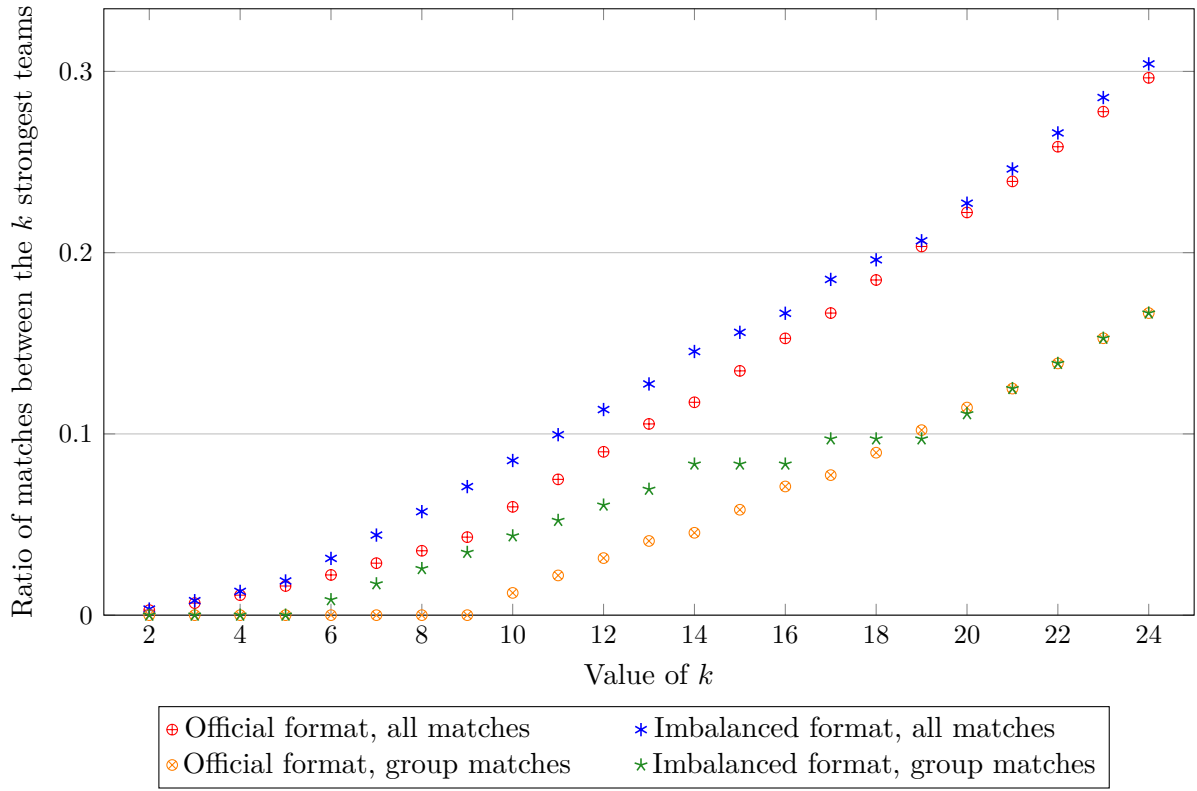


This finding is supported by Figure~\ref{Fig6}, which shows the share of matches between the $k$ strongest teams up to $k=24$. The imbalanced format clearly outperforms the official. Even the expected number of matches (not only their proportion) is higher for all values of $k$ up to 18. Thus, the novel design contains more matches between the top teams, contrary to the lower number of matches.

\begin{figure}[t!]
\centering

\begin{tikzpicture}
\begin{axis}[
name = axis1,
xlabel = Elo rating of the national team,
x label style = {font=\small},
ylabel = Probability of stakeless match if the team qualified,
y label style = {font=\small},
width = \textwidth,
height = 0.6\textwidth,
ymajorgrids = true,
ymin = 0,
ymax = 1.02,
max space between ticks = 50,
legend style = {font=\small,at={(0,-0.15)},anchor=north west,legend columns=2},
legend entries = {{Official format, two teams qualify ($s_i^{A,\min}$)}$\quad \: \: \:$, {Imbalanced format, Tier 1 groups}, {Official format, three teams qualify ($s_i^{A,\max}$)}$\quad$,{Imbalanced format, Tier 2 groups}}
]
\addplot [red, mark=oplus, only marks, mark size=2pt, mark options={solid,thin}] coordinates {
(1892,0.316209)
(1870,0.295098)
(1827,0.257564)
(2157,0.580056)
(2149,0.575332)
(2061,0.489656)
(2015,0.436073)
(1997,0.420577)
(1996,0.413513)
(1975,0.391354)
(1960,0.37578)
(1957,0.372757)
(1956,0.324866)
(1943,0.315212)
(1923,0.296445)
(1912,0.286343)
(1873,0.237048)
(1871,0.246209)
(1863,0.243138)
(1855,0.233746)
(1853,0.232059)
(1819,0.189303)
(1812,0.198347)
(1785,0.174319)
(1784,0.136133)
(1781,0.140317)
(1780,0.12677)
(1771,0.133218)
(1754,0.110397)
(1743,0.105679)
(1731,0.088707)
(1723,0.091783)
(1690,0.072516)
(1681,0.067844)
(1671,0.062605)
(1669,0.062006)
(1662,0.051695)
(1658,0.050831)
(1635,0.040598)
(1629,0.039464)
(1595,0.027637)
(1594,0.02763)
(1592,0.027425)
(1588,0.02722)
(1585,0.025907)
(1570,0.02389)
(1678.753,0.064774)
(1549.413,0.019902)
};
\addplot [blue, mark=asterisk, only marks, mark size=2.5pt, mark options={solid,semithick}] coordinates {
(1892,0.107601)
(1870,0.097875)
(1827,0.080988)
(2157,0.264286)
(2149,0.258723)
(2061,0.194884)
(2015,0.165632)
(1997,0.154849)
(1996,0.122678)
(1975,0.112501)
(1960,0.10575)
(1957,0.104386)
(1956,0.097534)
(1943,0.099374)
(1923,0.09153)
(1912,0.088709)
(1784,0.026057)
(1781,0.025291)
(1780,0.022202)
(1771,0.026794)
(1754,0.019018)
(1743,0.017626)
(1731,0.01191)
(1723,0.014547)
(1662,0.005574)
(1658,0.008841)
(1635,0.007688)
(1629,0.006797)
(1595,0.002795)
(1594,0.004956)
(1592,0.00504)
(1588,0.0046)
};
\addplot [orange, mark=otimes, only marks, mark size=2pt, mark options={solid,thin}] coordinates {
(1892,0.631353)
(1870,0.603044)
(1827,0.546982)
(2157,0.910643)
(2149,0.909254)
(2061,0.830925)
(2015,0.777121)
(1997,0.758766)
(1996,0.749576)
(1975,0.725978)
(1960,0.706357)
(1957,0.701714)
(1956,0.652424)
(1943,0.637598)
(1923,0.612383)
(1912,0.597204)
(1873,0.529526)
(1871,0.539029)
(1863,0.528252)
(1855,0.515574)
(1853,0.514276)
(1819,0.454091)
(1812,0.455638)
(1785,0.416367)
(1784,0.400047)
(1781,0.399122)
(1780,0.381233)
(1771,0.387043)
(1754,0.347837)
(1743,0.34305)
(1731,0.314221)
(1723,0.31825)
(1690,0.27879)
(1681,0.259584)
(1671,0.250296)
(1669,0.245768)
(1662,0.171522)
(1658,0.175124)
(1635,0.145733)
(1629,0.149342)
(1595,0.117014)
(1594,0.113559)
(1592,0.112704)
(1588,0.119071)
(1585,0.107683)
(1570,0.105368)
(1678.753,0.200631)
(1549.413,0.094252)
}; 
\addplot [ForestGreen, mark=triangle, only marks, mark size=2pt, mark options={solid,thick}] coordinates {
(1873,0.14097)
(1871,0.137792)
(1863,0.13399)
(1855,0.129098)
(1853,0.109217)
(1819,0.09551)
(1812,0.090564)
(1785,0.079117)
(1690,0.031173)
(1681,0.028268)
(1671,0.026157)
(1669,0.025885)
(1585,0.012601)
(1570,0.010893)
(1678.753,0.028855)
(1549.413,0.009569)
};
\end{axis}
\end{tikzpicture}

\captionsetup{justification=centerfirst}
\caption{The probability of a stakeless match played by an already qualified team in the two formats for the 2026 FIFA World Cup \\ \vspace{0.25cm}
\footnotesize{\emph{Notes}: The two winners of the play-offs are represented by their expected Elo ratings. \\
The Elo ratings of the three hosts are increased by 100 due to home advantage.}}
\label{Fig7}
\end{figure}


Finally, Figure~\ref{Fig7} conveys probably the most important message of our paper.
It provides the probability of a stakeless match played by a team that has already qualified either for the Round of 32 in the official format, or for the Round of 16 in Tier 1 groups of the imbalanced format, or for the knockout round play-offs in Tier 2 groups of the imbalanced format. In the proposed format, this value remains robustly below the lower bound $S_i^{A,\min}$ in the official format for any team.
Crucially, the probability of a stakeless game played after qualifying for the knockout stage is reduced by at least 17 percentage points (compared to the lower bound in the official format) for the 16 strongest teams, including the three host nations.

On the other hand, the probability of a stakeless game in the official 2026 FIFA World Cup is at least 28\% and can reach 60\% for the 16th strongest team. Furthermore, the gradually increasing lower and upper bounds exceed 58\% and 91\%, respectively, for the best team, Argentina.
Consequently, our imbalanced design would be highly efficient in avoiding the seriously non-competitive matches such as the ones played by three strong national teams (Brazil, France, Portugal) in the 2022 FIFA World Cup---all of them lost by the favourite.

\section{Conclusions} \label{Sec5}

The design of the 2026 FIFA World Cup has seen a significant reform with the expansion to 48 teams. The original plan, with 16 groups of three teams each, received serious criticism \citep{Guyon2020a, Stronka2024}, which prompted FIFA to choose a new format in 2023. The revised design contains 12 groups of four teams each, followed by a knockout stage starting from the Round of 32. Although the risk of collusion is decreased and at least three matches are guaranteed for each team, the qualification of two-thirds of the teams and the participation of 16 additional weak teams imply that non-competitive or stakeless matches are highly likely to occur, especially for the strongest teams, as demonstrated by our simulations (see Table~\ref{Table3} and Figure~\ref{Fig7}).

Inspired by examples from handball, ice hockey, and water polo, we have suggested an alternative design based on imbalanced groups. The basic idea is to divide the 48 teams into two divisions: eight groups of stronger teams (Tier 1) and four groups of weaker teams (Tier 2). While Tier 1 group winners directly qualify for the Round of 16, Tier 2 group winners and runners-up compete in a play-off round against Tier 1 runners-up.
Future studies should investigate whether this design can improve efficacy, the accuracy of the ranking since existing works \citep{LasekGagolewski2018, SziklaiBiroCsato2022} do not consider imbalanced groups.

The imbalanced group format offers several advantages over the official design.
First, it substantially reduces the proportion of stakeless matches, including the more costly cases.
Second, it contains fewer matches, especially for the strongest teams whose players have the highest workload at the end of the season.
Third, it ensures that the top teams face stronger opponents on average.

The current study demonstrates how the format of the 2026 FIFA World Cup can be improved to contain more competitive matches between top teams. The proposed design with imbalanced groups offers a viable alternative that maintains fairness while maximises attractiveness and excitement. These insights contribute to the rapidly growing literature on tournament design and can inform future discussions on optimising competitions where the strengths of the teams vary to a great extent.

\section*{Acknowledgements}

Three anonymous reviewers provided useful comments on an earlier draft. \\
The research was supported by the National Research, Development and Innovation Office under Grants Advanced 152220 and FK 145838, and the J\'anos Bolyai Research Scholarship of the Hungarian Academy of Sciences.

\bibliographystyle{apalike} 
\bibliography{All_references}

\begin{thebibliography}{}

\bibitem[Appleton, 1995]{Appleton1995}
Appleton, D.~R. (1995).
\newblock May the best man win?
\newblock {\em Journal of the Royal Statistical Society: Series D (The
  Statistician)}, 44(4):529--538.

\bibitem[Bengtsson et~al., 2013]{BengtssonEkstrandHagglund2013}
Bengtsson, H., Ekstrand, J., and H{\"a}gglund, M. (2013).
\newblock Muscle injury rates in professional football increase with fixture
  congestion: an 11-year follow-up of the {UEFA} {C}hampions {L}eague injury
  study.
\newblock {\em British Journal of Sports Medicine}, 47(12):743--747.

\bibitem[Brandes et~al., 2025]{BrandesMarmullaSmokovic2025}
Brandes, U., Marmulla, G., and Smokovic, I. (2025).
\newblock Efficient computation of tournament winning probabilities.
\newblock {\em Journal of Sports Analytics}, 11.
\newblock {DOI}:
  \href{https://doi.org/10.1177/22150218251313905}{10.1177/22150218251313905}.

\bibitem[Cea et~al., 2020]{CeaDuranGuajardoSureSiebertZamorano2020}
Cea, S., Dur{\'a}n, G., Guajardo, M., Saur{\'e}, D., Siebert, J., and Zamorano,
  G. (2020).
\newblock An analytics approach to the {FIFA} ranking procedure and the {W}orld
  {C}up final draw.
\newblock {\em Annals of Operations Research}, 286(1-2):119--146.

\bibitem[Chater et~al., 2021]{ChaterArrondelGayantLaslier2021}
Chater, M., Arrondel, L., Gayant, J.-P., and Laslier, J.-F. (2021).
\newblock Fixing match-fixing: Optimal schedules to promote competitiveness.
\newblock {\em European Journal of Operational Research}, 294(2):673--683.

\bibitem[Csat{\'o}, 2020]{Csato2020b}
Csat{\'o}, L. (2020).
\newblock Optimal tournament design: {L}essons from the men's handball
  {C}hampions {L}eague.
\newblock {\em Journal of Sports Economics}, 21(8):848--868.

\bibitem[Csat{\'o}, 2021]{Csato2021a}
Csat{\'o}, L. (2021).
\newblock {\em Tournament Design: How Operations Research Can Improve Sports
  Rules}.
\newblock Palgrave Pivots in Sports Economics. Palgrave Macmillan, Cham,
  Switzerland.

\bibitem[Csat{\'o}, 2022a]{Csato2022d}
Csat{\'o}, L. (2022a).
\newblock The effects of draw restrictions on knockout tournaments.
\newblock {\em Journal of Quantitative Analysis in Sports}, 18(4):227--239.

\bibitem[Csat{\'o}, 2022b]{Csato2022a}
Csat{\'o}, L. (2022b).
\newblock Quantifying incentive (in)compatibility: {A} case study from sports.
\newblock {\em European Journal of Operational Research}, 302(2):717--726.

\bibitem[Csat{\'o}, 2023a]{Csato2023d}
Csat{\'o}, L. (2023a).
\newblock Group draw with unknown qualified teams: {A} lesson from the 2022
  {FIFA} {W}orld {C}up.
\newblock {\em International Journal of Sports Science \& Coaching},
  18(2):539--551.

\bibitem[Csat{\'o}, 2023b]{Csato2023a}
Csat{\'o}, L. (2023b).
\newblock How to avoid uncompetitive games? {T}he importance of tie-breaking
  rules.
\newblock {\em European Journal of Operational Research}, 307(3):1260--1269.

\bibitem[Csat{\'o}, 2023c]{Csato2023c}
Csat{\'o}, L. (2023c).
\newblock Quantifying the unfairness of the 2018 {FIFA} {W}orld {C}up
  qualification.
\newblock {\em International Journal of Sports Science \& Coaching},
  18(1):183--196.

\bibitem[Csat{\'o}, 2025a]{Csato2025c}
Csat{\'o}, L. (2025a).
\newblock The fairness of the group draw for the {FIFA} {W}orld {C}up.
\newblock {\em International Journal of Sports Science \& Coaching},
  20(2):554--567.

\bibitem[Csat{\'o}, 2025b]{Csato2025d}
Csat{\'o}, L. (2025b).
\newblock On head-to-head results as tie-breaker and consequent opportunities
  for collusion.
\newblock {\em IMA Journal of Management Mathematics}, 36(2):215--230.

\bibitem[Csat{\'o} et~al., 2025]{CsatoKissSzadoczki2025}
Csat{\'o}, L., Kiss, L.~M., and Sz{\'a}doczki, {\relax Zs}. (2025).
\newblock The allocation of {FIFA} {W}orld {C}up slots based on the ranking of
  confederations.
\newblock {\em Annals of Operations Research}, 344(1):153--173.

\bibitem[Csat{\'o} et~al., 2024]{CsatoMolontayPinter2024}
Csat{\'o}, L., Molontay, R., and Pint{\'e}r, J. (2024).
\newblock Tournament schedules and incentives in a double round-robin
  tournament with four teams.
\newblock {\em International Transactions in Operational Research},
  31(3):1486--1514.

\bibitem[Dellal et~al., 2015]{DellalLago-PenasReyChamariOrhant2015}
Dellal, A., Lago-Pe{\~n}as, C., Rey, E., Chamari, K., and Orhant, E. (2015).
\newblock The effects of a congested fixture period on physical performance,
  technical activity and injury rate during matches in a professional soccer
  team.
\newblock {\em British Journal of Sports Medicine}, 49(6):390--394.

\bibitem[Devriesere et~al., 2025a]{DevriesereCsatoGoossens2025}
Devriesere, K., Csat{\'o}, L., and Goossens, D. (2025a).
\newblock Tournament design: A review from an operational research perspective.
\newblock {\em European Journal of Operational Research}, 324(1):1--21.

\bibitem[Devriesere et~al., 2025b]{DevriesereGoossensSpieksma2025}
Devriesere, K., Goossens, D., and Spieksma, F. (2025b).
\newblock Evaluating competitiveness in {UEFA}'s new {C}hampions {L}eague
  format.
\newblock Manuscript. {DOI}:
  \href{https://doi.org/10.48550/arXiv.2508.08290}{10.48550/arXiv.2508.08290}.

\bibitem[EHF, 2014]{EHF2014a}
EHF (2014).
\newblock The 119th meeting of the {EHF} {E}xecutive {C}ommittee.
\newblock 21 March.
  \url{http://www.eurohandball.com/article/018817/EHF+Executive+meets+in+Vienna}.

\bibitem[Ekstrand et~al., 2004]{EkstrandWaldenHagglund2004}
Ekstrand, J., Wald{\'e}n, M., and H{\"a}gglund, M. (2004).
\newblock A congested football calendar and the wellbeing of players:
  correlation between match exposure of {E}uropean footballers before the
  {W}orld {C}up 2002 and their injuries and performances during that {W}orld
  {C}up.
\newblock {\em British Journal of Sports Medicine}, 38(4):493--497.

\bibitem[FIFA, 2017a]{FIFA2017c}
FIFA (2017a).
\newblock Close-up on {F}inal {D}raw procedures.
\newblock 27 November.
  \url{https://web.archive.org/web/20171127150059/http://www.fifa.com/worldcup/news/y=2017/m=11/news=close-up-on-final-draw-procedures-2921440.html}.

\bibitem[FIFA, 2017b]{FIFA2017d}
FIFA (2017b).
\newblock Unanimous decision expands {FIFA} {W}orld {C}up\textsuperscript{{TM}}
  to 48 teams from 2026.
\newblock 10 January.
  \url{https://web.archive.org/web/20170110231324/http://www.fifa.com/about-fifa/news/y=2017/m=1/news=fifa-council-unanimously-decides-on-expansion-of-the-fifa-world-cuptm--2863100.html}.

\bibitem[FIFA, 2018]{FIFA2018c}
FIFA (2018).
\newblock Revision of the {FIFA} / {C}oca-{C}ola {W}orld {R}anking.
\newblock
  \url{https://digitalhub.fifa.com/m/f99da4f73212220/original/edbm045h0udbwkqew35a-pdf.pdf}.

\bibitem[FIFA, 2022a]{FIFA2022a}
FIFA (2022a).
\newblock {\em {D}raw procedures. {FIFA} {W}orld {C}up {Q}atar
  2022\textsuperscript{{TM}}}.
\newblock
  \url{https://digitalhub.fifa.com/m/2ef762dcf5f577c6/original/Portrait-Master-Template.pdf}.

\bibitem[FIFA, 2022b]{FIFA2022d}
FIFA (2022b).
\newblock {\em {FIFA} {S}tatues}.
\newblock May.
  \url{https://digitalhub.fifa.com/m/3815fa68bd9f4ad8/original/FIFA_Statutes_2022-EN.pdf}.

\bibitem[FIFA, 2023]{FIFA2023b}
FIFA (2023).
\newblock {FIFA} {C}ouncil approves international match calendars.
\newblock 14 March.
  \url{https://www.fifa.com/about-fifa/organisation/fifa-council/media-releases/fifa-council-approves-international-match-calendars}.

\bibitem[FIFA, 2025]{FIFA2025}
FIFA (2025).
\newblock {\em Regulations for the {FIFA} {W}orld {C}up
  26\textsuperscript{TM}}.
\newblock
  \url{https://digitalhub.fifa.com/m/636f5c9c6f29771f/original/FWC2026_regulations_EN.pdf}.

\bibitem[Folgado et~al., 2015]{FolgadoDuarteMarquesSampaio2015}
Folgado, H., Duarte, R., Marques, P., and Sampaio, J. (2015).
\newblock The effects of congested fixtures period on tactical and physical
  performance in elite football.
\newblock {\em Journal of Sports Sciences}, 33(12):1238--1247.

\bibitem[{Football rankings}, 2020]{FootballRankings2020}
{Football rankings} (2020).
\newblock Simulation of scheduled football matches.
\newblock 28 December.
  \url{http://www.football-rankings.info/2020/12/simulation-of-scheduled-football-matches.html}.

\bibitem[G{\'a}squez and Royuela, 2016]{GasquezRoyuela2016}
G{\'a}squez, R. and Royuela, V. (2016).
\newblock The determinants of international football success: A panel data
  analysis of the {E}lo rating.
\newblock {\em Social Science Quarterly}, 97(2):125--141.

\bibitem[Gerchak et~al., 1995]{GerchakMausserMagazine1995}
Gerchak, Y., Mausser, H.~E., and Magazine, M.~J. (1995).
\newblock The evolution of draft lotteries in professional sports: Back to
  moral hazard?
\newblock {\em Interfaces}, 25(6):30--38.

\bibitem[Guajardo and Krumer, 2024]{GuajardoKrumer2024}
Guajardo, M. and Krumer, A. (2024).
\newblock Tournament design for a {FIFA} {W}orld {C}up with 12 four-team
  groups: Every win matters.
\newblock In Breuer, M. and Forrest, D., editors, {\em The Palgrave Handbook on
  the Economics of Manipulation in Sport}, pages 207--230. Palgrave Macmillan,
  Cham, Switzerland.

\bibitem[Guyon, 2015]{Guyon2015a}
Guyon, J. (2015).
\newblock Rethinking the {FIFA} {W}orld {C}up\textsuperscript{{TM}} final draw.
\newblock {\em Journal of Quantitative Analysis in Sports}, 11(3):169--182.

\bibitem[Guyon, 2018]{Guyon2018a}
Guyon, J. (2018).
\newblock What a fairer 24 team {UEFA} {E}uro could look like.
\newblock {\em Journal of Sports Analytics}, 4(4):297--317.

\bibitem[Guyon, 2020]{Guyon2020a}
Guyon, J. (2020).
\newblock Risk of collusion: {W}ill groups of 3 ruin the {FIFA} {W}orld {C}up?
\newblock {\em Journal of Sports Analytics}, 6(4):259--279.

\bibitem[Guyon, 2022]{Guyon2022a}
Guyon, J. (2022).
\newblock ``{C}hoose your opponent'': A new knockout design for hybrid
  tournaments.
\newblock {\em Journal of Sports Analytics}, 8(1):9--29.

\bibitem[Gyimesi, 2024]{Gyimesi2024}
Gyimesi, A. (2024).
\newblock Competitive balance in the post-2024 {C}hampions {L}eague and the
  {E}uropean {S}uper {L}eague: {A} simulation study.
\newblock {\em Journal of Sports Economics}, 25(6):707--734.

\bibitem[IIHF, 2024]{IIHF2024}
IIHF (2024).
\newblock 2025 {W}omen's {W}orld {C}hampionship {T}ournament {I}nfo.
\newblock 8 October.
  \url{https://www.iihf.com/en/events/2025/ww/tournamentinfo/62150/tournament_info}.

\bibitem[Kendall and Lenten, 2017]{KendallLenten2017}
Kendall, G. and Lenten, L.~J.~A. (2017).
\newblock When sports rules go awry.
\newblock {\em European Journal of Operational Research}, 257(2):377--394.

\bibitem[Krumer and Lechner, 2017]{KrumerLechner2017}
Krumer, A. and Lechner, M. (2017).
\newblock First in first win: Evidence on schedule effects in round-robin
  tournaments in mega-events.
\newblock {\em European Economic Review}, 100:412--427.

\bibitem[Krumer and Moreno-Ternero, 2023]{KrumerMoreno-Ternero2023}
Krumer, A. and Moreno-Ternero, J. (2023).
\newblock The allocation of additional slots for the {FIFA} {W}orld {C}up.
\newblock {\em Journal of Sports Economics}, 24(7):831--850.

\bibitem[Laica et~al., 2021]{LaicaLauberSahm2021}
Laica, C., Lauber, A., and Sahm, M. (2021).
\newblock Sequential round-robin tournaments with multiple prizes.
\newblock {\em Games and Economic Behavior}, 129:421--448.

\bibitem[Laliena and L{\'o}pez, 2019]{LalienaLopez2019}
Laliena, P. and L{\'o}pez, F.~J. (2019).
\newblock Fair draws for group rounds in sport tournaments.
\newblock {\em International Transactions in Operational Research},
  26(2):439--457.

\bibitem[Laliena and L{\'o}pez, 2025]{LalienaLopez2025}
Laliena, P. and L{\'o}pez, F.~J. (2025).
\newblock Draw procedures for balanced 3-team group rounds in sports
  competitions.
\newblock {\em Annals of Operations Research}, 346(3):2065--2092.

\bibitem[Lapr{\' e} and Palazzolo, 2023]{LaprePalazzolo2023}
Lapr{\' e}, M.~A. and Palazzolo, E.~M. (2023).
\newblock The evolution of seeding systems and the impact of imbalanced groups
  in {FIFA} {M}en's {W}orld {C}up tournaments 1954--2022.
\newblock {\em Journal of Quantitative Analysis in Sports}, 19(4):317--332.

\bibitem[Lasek and Gagolewski, 2018]{LasekGagolewski2018}
Lasek, J. and Gagolewski, M. (2018).
\newblock The efficacy of league formats in ranking teams.
\newblock {\em Statistical Modelling}, 18(5-6):411--435.

\bibitem[Lasek et~al., 2013]{LasekSzlavikBhulai2013}
Lasek, J., Szl{\'a}vik, Z., and Bhulai, S. (2013).
\newblock The predictive power of ranking systems in association football.
\newblock {\em International Journal of Applied Pattern Recognition},
  1(1):27--46.

\bibitem[{Le Monde}, 2022]{LeMonde2022}
{Le Monde} (2022).
\newblock {W}orld {C}up 2022: Dramatic 1-0 victory for {T}unisia over {F}rance,
  but {L}es {B}leus still win group.
\newblock 30 November.
  \url{https://www.lemonde.fr/en/football/article/2022/11/30/world-cup-2022-dramatic-1-0-victory-for-tunisia-over-france-but-les-bleus-still-win-group_6006217_130.html}.

\bibitem[LEN, 2023]{LEN2023}
LEN (2023).
\newblock 2024 {E}uropean {W}ater {P}olo {C}hampionships draw: Everything you
  need to know.
\newblock 8 September.
  \url{https://www.len.eu/2024-european-water-polo-championships-draw-everything-you-need-to-know/}.

\bibitem[L{\'o}pez-Valenciano et~al.,
  2020]{Lopez-ValencianoRuiz-PerezGarcia-GomezVera-GarciaCroixMyerAyala2020}
L{\'o}pez-Valenciano, A., Ruiz-P{\'e}rez, I., Garcia-G{\'o}mez, A.,
  Vera-Garcia, F.~J., Croix, M.~D.~S., Myer, G.~D., and Ayala, F. (2020).
\newblock Epidemiology of injuries in professional football: a systematic
  review and meta-analysis.
\newblock {\em British Journal of Sports Medicine}, 54(12):711--718.

\bibitem[Maher, 1982]{Maher1982}
Maher, M.~J. (1982).
\newblock Modelling association football scores.
\newblock {\em Statistica Neerlandica}, 36(3):109--118.

\bibitem[Medcalfe, 2024]{Medcalfe2024}
Medcalfe, S. (2024).
\newblock Behavioral responses to sporting contest design: A review of the
  literature.
\newblock {\em American Behavioral Scientist}, in press.
\newblock {DOI}:
  \href{https://doi.org/10.1177/00027642241235814}{10.1177/00027642241235814}.

\bibitem[Palacios-Huerta, 2023]{Palacios-Huerta2023b}
Palacios-Huerta, I. (2023).
\newblock The beautiful dataset.
\newblock Manuscript. {DOI}:
  \href{https://doi.org/10.2139/ssrn.4665889}{10.2139/ssrn.4665889}.

\bibitem[Renn{\'o}-Costa, 2023]{Renno-Costa2023}
Renn{\'o}-Costa, C. (2023).
\newblock A double-elimination format for a 48-team {FIFA} {W}orld {C}up.
\newblock Manuscript. {DOI}:
  \href{https://doi.org/10.48550/arXiv.2301.03411}{10.48550/arXiv.2301.03411}.

\bibitem[Roberts and Rosenthal, 2024]{RobertsRosenthal2024}
Roberts, G.~O. and Rosenthal, J.~S. (2024).
\newblock Football group draw probabilities and corrections.
\newblock {\em The Canadian Journal of Statistics}, 52(3):659--677.

\bibitem[Scelles et~al., 2024]{ScellesFrancoisValenti2024}
Scelles, N., Fran{\c{c}}ois, A., and Valenti, M. (2024).
\newblock Impact of the {UEFA} {N}ations {L}eague on competitive balance,
  competitive intensity, and fairness in {E}uropean men's national team
  football.
\newblock {\em International Journal of Sport Policy and Politics},
  16(3):519--537.

\bibitem[Schmidt, 2024]{Schmidt2024}
Schmidt, M.~B. (2024).
\newblock On the incentive structure of tournaments: Evidence from the
  {N}ational {B}asketball {A}ssociation's draft lottery.
\newblock {\em Eastern Economic Journal}, 50(3):400--429.

\bibitem[Stronka, 2024]{Stronka2024}
Stronka, W. (2024).
\newblock Demonstration of the collusion risk mitigation effect of random
  tie-breaking and dynamic scheduling.
\newblock {\em Sports Economics Review}, 5:100025.

\bibitem[Szczecinski and Roatis, 2022]{SzczecinskiRoatis2022}
Szczecinski, L. and Roatis, I.-I. (2022).
\newblock {FIFA} ranking: Evaluation and path forward.
\newblock {\em Journal of Sports Analytics}, 8(4):231--250.

\bibitem[Sziklai et~al., 2022]{SziklaiBiroCsato2022}
Sziklai, B.~R., Bir\'o, P., and Csat{\'o}, L. (2022).
\newblock The efficacy of tournament designs.
\newblock {\em Computers \& Operations Research}, 144:105821.

\bibitem[Szymanski, 2003]{Szymanski2003}
Szymanski, S. (2003).
\newblock The economic design of sporting contests.
\newblock {\em Journal of Economic Literature}, 41(4):1137--1187.

\bibitem[Truta, 2018]{Truta2018}
Truta, T.~M. (2018).
\newblock {FIFA} does it right: 2026 {FIFA} {W}orld {C}up does not increase the
  number of non-competitive matches.
\newblock Manuscript. DOI:
  \href{https://doi.org/10.48550/arXiv.1808.05858}{10.48550/arXiv.1808.05858}.

\bibitem[van Eetvelde and Ley, 2019]{vanEetveldeLey2019}
van Eetvelde, H. and Ley, C. (2019).
\newblock Ranking methods in soccer.
\newblock In Kenett, R.~S., Longford, T.~N., Piegorsch, W., and Ruggeri, F.,
  editors, {\em Wiley StatsRef: Statistics Reference Online}, pages 1--9.
  Springer, Hoboken, New Jersey, USA.

\bibitem[Wright, 2009]{Wright2009}
Wright, M. (2009).
\newblock 50 years of {OR} in sport.
\newblock {\em Journal of the Operational Research Society}, 60(Supplement
  1):S161--S168.

\bibitem[Wright, 2014]{Wright2014}
Wright, M. (2014).
\newblock {OR} analysis of sporting rules -- {A} survey.
\newblock {\em European Journal of Operational Research}, 232(1):1--8.

\end{thebibliography}

\end{document}